\documentclass[letter,onecolumn,12pt]{IEEEtran2}

\usepackage{amssymb}
\usepackage{amsmath}
\usepackage{amsthm}
\usepackage{latexsym}
\usepackage{epsfig}
\usepackage{graphicx,subfigure}
\usepackage{color}
\usepackage{float}
\usepackage{cite}
\usepackage{enumerate}

\long\def\comment#1{}

\renewcommand{\thesection}{\arabic{section}}
\renewcommand{\thesubsection}{\arabic{section}.\arabic{subsection}}

\newtheorem{theorem}{Theorem}[section]
\newtheorem{remark}{Remark}[section]
\newtheorem{lemma}{Lemma}[section]
\newtheorem{example}{Example}[section]

\newtheorem{corollary}{Corollary}[section]
\newtheorem{proposition}{Proposition}[section]
\newtheorem{definition}{Definition}[section]

\newtheorem{thmsub}{Theorem}[subsection]
\newtheorem{remsub}{Remark}[subsection]
\newtheorem{lemsub}{Lemma}[subsection]
\newtheorem{exmsub}{Example}[subsection]

\newtheorem{corsub}{Corollary}[subsection]
\newtheorem{propsub}{Proposition}[subsection]
\newtheorem{defsub}{Definition}[subsection]

\begin{document}

\title{Detectability, Observability and Lyapunov-Type Theorems of Linear Discrete Time-Varying Stochastic
Systems with Multiplicative Noise
\thanks{This work was supported in part by NSF of China under Grant 61573227 and a research grant from the Australian Research Council.}
}

\author{Weihai Zhang$^\dag$,~
        Wei Xing~Zheng, ~and~ Bor-Sen Chen
\thanks{W. Zhang is with the College of Electrical Engineering and Automation,
Shandong University of Science and Technology, Qingdao 266590, China.}
\thanks{W. X. Zheng is with the School of Computing, Engineering and Mathematics,
Western Sydney University, Sydney, NSW 2751, Australia.}
\thanks{B. S. Chen is with the Department of Electrical Engineering,
National Tsing Hua University, Hsinchu 30013, Taiwan.}
\thanks{$\dag$ Corresponding author. Email: w\_hzhang@163.com.}
}

%\markboth{Submitted to IEEE Transactions on Automatic Control}%
%{Shell \MakeLowercase{\textit{et al.}}: Bare Demo of IEEEtran.cls for Journals}

\maketitle

\begin{abstract}
The objective of this paper is to study detectability, observability
and related Lyapunov-type theorems of linear discrete-time time-varying
stochastic systems with multiplicative noise. Some new concepts such
as uniform detectability, ${\cal K}^{\infty}$-exact detectability
(resp. ${\cal K}^{WFT}$-exact detectability, ${\cal K}^{FT}$-exact
detectability, ${\cal K}^{N}$-exact detectability) and ${\cal
K}^{\infty}$-exact observability (resp. ${\cal K}^{WFT}$-exact
observability, ${\cal K}^{FT}$-exact observability, ${\cal
K}^{N}$-exact observability) are introduced, respectively, and nice
properties associated with uniform detectability, exact
detectability and exact observability are also obtained. Moreover,
some Lyapunov-type theorems associated with generalized Lyapunov
equations and exponential stability in mean square sense are
presented under uniform detectability, ${\cal K}^{N}$-exact
observability and ${\cal K}^{N}$-exact detectability, respectively.
\end{abstract}


%{\bf AMS subject classifications.} 47N70, 15A24, 93B07, 93E15, 93D15, 93C05

%%%%%%%%%%%%%%%%%%%%%%%%%%%%%%%%%%%%%%%%%%%%%%%%%%%%%%%%%%%%%%%%%%%%%%%%%%%%%%%%
\section{Introduction}
%%%%%%%%%%%%%%%%%%%%%%%%%%%%%%%%%%%%%%%%%%%%%%%%%%%%%%%%%%%%%%%%%%%%%%%%%%%%%%%%

It is well-known that observability and detectability are
fundamental concepts in system analysis and synthesis; see, e.g.,
\cite{bdo, c2ddss, c1ss, rugh, zhang1, c3dd, seo, Damm1, shenl,
lguo,lizao, Halanay, liwang}. In the linear system theory,
detectability is a weaker concept than observability, since it
describes the fact that all unobservable states are asymptotically
stable. Over the last two decades, the classical detectability in
the linear system theory has been extended to stochastic systems in
different ways. For example, the definition of stochastic
detectability for time-invariant It\^o stochastic systems can be
found in \cite{damm2004,c2ddss}, which is dual to mean square
stabilization. In \cite{xshf, c3dd,Damm1}, the notions of exact
observability and exact detectability were presented for It\^o
stochastic systems, which led to the stochastic
Popov-Belevitch-Hautus (PBH) criteria like those for deterministic
systems. Another natural concept of detectability for It\^o
stochastic systems was given in \cite{Damm1} based on the idea that
any non-observed states corresponded to stable models of the system.
In \cite{liwang}, the exact detectability in \cite{c3dd} and
detectability in \cite{Damm1} were proved to be equivalent, and a
unified treatment was proposed for detectability and observability
of It\^o stochastic systems. Based on the standard notions of
detectability and observability for time-varying linear systems
\cite{bdo,peters}, studied in \cite{lizao} were detectability and
observability of discrete time-invariant stochastic systems as well
as the properties of Lyapunov equations. Recently, the exact
detectability and observability were extended to stochastic systems
with Markov jumps and multiplicative noise in \cite{ni,deng1,
shenl,zhangtan}.

As it is well known that the classical Lyapunov theorem is very
essential in stability theory, which asserts that if a matrix $F$
is Schur stable, then for any $Q\ge 0$, the classical Lyapunov
equation $-P+F^TPF+Q=0$ admits a unique solution $P\ge 0$;
Conversely, if $(F,Q)$ is detectable, $Q\ge 0$, and the Lyapunov
equation $-P+F^TPF+Q=0$ admits a unique solution $P\ge 0$, then $F$
is Schur stable. The classical Lyapunov theorem was generalized to
deterministic time-varying systems in \cite{bdo} and will be
extended to stochastic time-varying systems in this paper under any
one assumption of uniform detectability, ${\cal K}^N$-exact
detectability and ${\cal K}^N$-exact observability.

Recently it has become known that discrete-time stochastic systems
with multiplicative noise are ideal models in the fields of investment
portfolio optimization \cite{dom}, system biology \cite{chens} and so on.
So the discrete-time stochastic $H_2/H_\infty$ control and filtering design
have been extensively studied in recent years; see, e.g.,
\cite{c1,zhang1, c2ddss, hzhang, dragan2013, vfcv} and the
references therein. As it is well-known, time-varying systems may be
utilized to model more realistic systems and are more challenging in
mathematics than time-invariant ones. So far, the majority of the
existing results is focused on detectability of time-invariant
systems only, except for a few about time-varying systems; see
\cite{bdo,Ungureanu1,Ungureanu2,xsdsffdf, dragan2015,
lguo,dragan,dragan2013,seo,jwlee}. Because linear time-invariant
systems are not sufficient to describe many practical phenomena,
this motivates researchers to study time-varying systems. In the
classical work \cite{bdo}, uniform detectability of the
deterministic linear discrete-time time-varying (LDTV) system
\begin{equation}
%% x_{k+1}=F_kx_k, \ y_k=H_kx_k \label{eq vfhv}
\left\{
\begin{split}
x_{k+1}&=F_kx_k, \ x_0\in {\cal R}^n \\
y_k&=H_kx_k, \ k=0,1,2, \cdots
\end{split}
\right. \label{eq vfhv}
\end{equation}
was defined and discussed. By the duality of stochastic
stabilizability, another definition called ``stochastic
detectability" was introduced in \cite{c2ddss} for LDTV Markov
systems, which is not equivalent to uniform detectability in
time-varying case.

Mainly motivated by the preceding discussion and the authors' series
works \cite{xshf,zhang1,xsdsffdf,c3dd}, this paper will study
detectability, observability and Lyapunov-type equation related
to LDTV stochastic systems with multiplicative noise. Firstly, the
classical uniform detectability of \cite{bdo} for such systems is
extended, and some properties on uniform detectability are obtained.
By means of our Lemma~\ref{lem:2.1} given later, we obtain the
observability Gramian matrix ${\cal O}_{k+s,k}$ and the state
transition matrix $\phi_{l,k}$, which are deterministic matrices and
easy to be applied in practice. Specifically, we prove an important
theorem that uniform detectability preserves invariance under an
output feedback control law, which is expected to be useful in
stochastic $H_2/H_\infty$ control. As an application, under the
assumption of uniform detectability, Lyapunov-type theorems on
stochastic stability are also presented.

Secondly, we extend exact detectability of linear continuous-time
stochastic It\^o systems \cite{Damm1, xshf} to LDTV systems. We
introduce four concepts called ${\cal K}^{N}$-exact detectability,
${\cal K}^{FT}$-exact detectability, ${\cal K}^{WFT}$-exact
detectability and ${\cal K}^\infty$-exact detectability, and they in
turn become weaker in the sense that the former implies the latter
in a sequence. Although in linear time-invariant system
\begin{equation}
\left\{
\begin{split}
x_{k+1}&=Fx_k, \ x_0\in {\cal R}^n \\
y_k&=Hx_k, \ \ k=0,1,2, \cdots
\end{split}
\right. \label{eq vgthgh}
\end{equation}
these four concepts are equivalent with $N=n-1$, but  they are
different from the others in the time-varying case, which reveals
the essential difference between time-invariant and time-varying
systems. It is shown that uniform detectability implies ${\cal
K}^\infty$-exact detectability (see Lemma~\ref{lem:3.1.3}), and
stochastic detectability \cite{c2ddss} implies the above four types
of exact detectability (see Proposition~\ref{prop:3.1.1} and
Remark~\ref{rem:3.1.3}). It seems that there is no inclusion
relation among uniform detectability, ${\cal K}^{N}$-exact
detectability, ${\cal K}^{FT}$-exact detectability and ${\cal
K}^{WFT}$-exact detectability, although they can be unified in the
linear discrete time-invariant systems \cite{lizao}. Two important
Lyapunov-type theorems under ${\cal K}^{N}$-exact detectability for
periodic systems are obtained (see
Theorems~\ref{thm:5.3.1}--\ref{thm:5.3.2}), which reveal the
important relation between the exponential stability and the
existence of positive definite solutions of generalized Lyapunov
equations (GLEs).

Parallel to various definitions on exact detectability, we also
introduce ${\cal K}^{N}$-exact observability, ${\cal K}^{FT}$-exact
observability, ${\cal K}^{WFT}$-exact observability and ${\cal
K}^\infty$-exact observability, which are respectively stronger than
${\cal K}^{N}$-exact detectability, ${\cal K}^{FT}$-exact
detectability, ${\cal K}^{WFT}$-exact detectability and ${\cal
K}^\infty$-exact detectability. For the linear time-invariant system
(\ref{eq vgthgh}), ${\cal K}^{n-1}$-, ${\cal K}^{FT}$-, ${\cal
K}^{WFT}$- and ${\cal K}^\infty$-exact observability are equivalent,
but they are different definitions for the linear time-varying
system (\ref{eq vfhv}). We present a rank criterion for ${\cal
K}^\infty$- and a criterion for ${\cal K}^{N}$-exact observability
 based on the Gramian matrix ${\cal O}_{k+N,k}$. Finally,
under the assumption of ${\cal K}^{N}$-exact observability, a
Lyapunov-type theorem is derived from Theorem~\ref{thm:5.3.1}.

The rest of the paper is organized as follows. In
Section~\ref{sec:UD}, we define uniform detectability and discuss
its properties. Lyapunov-type theorems are given under uniform
detectability. Section~\ref{sec:EDEO} introduces some new concepts
about exact detectability and exposes nice properties. This section
also presents Lyapunov-type stability theorems based on ${\cal
K}^{N}$-exact detectability. Moreover, the relation among uniform
detectability, exact detectability and stochastic detectability is
clarified via some examples. Section~\ref{sec:TLTT} introduces
various definitions for exact observability, which are stronger
than those of Section-\ref{sec:EDEO-1}. Section~\ref{sec:LTTED}
provides some comments on this study. Finally,
Section~\ref{sec:Concl} concludes the paper with some remarks.

\textit{Notation:} \ ${\cal R}^n$: the set of all real
$n$-dimensional vectors. ${\cal S}_n$: the set of all $n\times n$
symmetric matrices whose entries may be complex numbers. ${\cal C}$:
the set of all complex numbers. $\mathcal{R}^{m\times n}$: the set
of all $m\times n$ real matrices. $\|x\|$: the norm of a vector or
matrix. $A>0$ (resp. $A\ge 0$): $A$ is a real symmetric positive
definite (resp. positive semi-definite) matrix. $I$: the identity
matrix. $\sigma(L)$: the spectrum set of the operator or matrix $L$.
$A^T$: the transpose of matrix $A$. ${\cal N}_{k_0}:=\{k_0,k_0+1,
k_0+2,\cdots, \}$, especially, ${\cal N}_1=\{1,2,\cdots, \}$, ${\cal
N}_0=\{0, 1,2,\cdots, \}$. $l^2_{{\cal F}_k}:=\{x(\omega): x \
\mbox{is} \ {\cal F}_k-\mbox{measurable}, E\|x\|^2<\infty\}$.

%%%%%%%%%%%%%%%%%%%%%%%%%%%%%%%%%%%%%%%%%%%%%%%%%%%%%%%%%%%%%%%%%%%%%%%%%%%%%%%%
\section{Uniform Detectability and Related Lyapunov-type Theorems\label{sec:UD}}
%%%%%%%%%%%%%%%%%%%%%%%%%%%%%%%%%%%%%%%%%%%%%%%%%%%%%%%%%%%%%%%%%%%%%%%%%%%%%%%%

In this section, we will define one important concept for LDTV
stochastic systems, called ``uniform detectability". And then, we
will obtain Lyapunov-type theorems under uniform detectability, which are
extensions of classical Lyapunov theorem.

\subsection{Uniform Detectability}\label{sec:UD-1}

Consider the following LDTV stochastic system
\begin{equation}
\left\{
\begin{split}
x_{k+1}&=F_kx_k+G_kx_kw_k, \ x_0\in {\cal R}^n\\
y_k&=H_kx_k, \ k=0,1,2,\cdots,
\end{split}
\right. \label{eq vcgtgh1}
\end{equation}
where $x_k$ is the $n$-dimensional state vector, $y_k$ is the
$m$-dimensional measurement output, $\{w_k\}_{k\ge 0}$ represents a
one-dimensional independent white noise process defined on the
filtered probability space $(\Omega, {\cal F}, {\cal F}_{k}, {\cal
P})$ with ${\cal F}_k=\sigma(w(0),\cdots, w(k))$. Assume that
$Ew_k=0$, $E[w_kw_j]=\delta_{kj}$, where $\delta_{kj}$ is a
Kronecker function defined by $\delta_{kj}=0$ for $k\ne j$ while
$\delta_{kj}=1$ for $k=j$. $x_0$ is assumed to be deterministic for
simplicity purposes, and $F_k$, $G_k$ and $H_k$ are time-varying
matrices of appropriate dimension. In practice, one is more
concerned about the $l^2_{{\cal F}_k}$-solution $\{x_k\}_{k\in {\cal
N}_0}$ of stochastic difference equation
\begin{equation}
x_{k+1}=F_kx_k+G_kx_kw_k, \ x_0\in {\cal R}^n . \label{eq vhvhc}
\end{equation}

\begin{definition}\label{def:2.1}\
The stochastic vector-valued sequence $\{{\tilde x}_k\}_{k\in {\cal
N}_0}$ is called a solution of system (\ref{eq vhvhc}) if (i)
${\tilde {x}}_0=x_0$; (ii) ${\tilde x}_k$ solves (\ref{eq vhvhc})
for $k=1,2,\cdots$; (iii) ${\tilde x}_k\in l^2_{{\cal F}_{k-1}}$,
where ${\cal F}_{-1}=\{\phi,\Omega\}$ is assumed to be a trivial
sigma algebra. System (\ref{eq vhvhc}) is said to have a unique
solution if for any two of its solutions $\{{\tilde x}_k\}_{k\in
{\cal N}_0}$ and $\{{\bar x}_k\}_{k\in {\cal N}_0}$, ${\cal
P}({\tilde x}_k={\bar x}_k, k\in {\cal N}_0)=1$.
\end{definition}

\begin{remark}\label{rem:2.1}\rm\
It can be found that, in most present literature, the condition (iii) in
Definition~\ref{def:2.1} is not particularly pointed out when defining
solutions of system (\ref{eq vhvhc}), which is in fact an essential requirement
as done in stochastic differential equations \cite{mao}. This makes an
fundamental difference of (\ref{eq vhvhc}) from deterministic difference
equations, as will be seen in the following examples.
\end{remark}

\begin{example}\label{exm:2.1}\rm\
It is easy to see that the following forward difference equation
$$
x_{k+1}=F_kx_k, \ x_0\in {\cal R}^n, \ k=0,1, \cdots, N
$$
always admits a unique solution on $[0,N+1]$. In addition, if
$F_k, k=0, 1,\cdots, N$ are nonsingular, then the backward difference
equation
$$
x_{k+1}=F_kx_k, \ x_{N+1}\in {\cal R}^n, \ k=0,1,\cdots, N
$$
also has a unique solution on $[0,N+1]$.
\end{example}

\begin{example}\label{exm:2.2}\rm\
Obviously, the linear stochastic difference equation (\ref{eq
vhvhc}) always has a unique $l^2_{{\cal F}_{k-1}}$-solution $x_k$ on
any interval $[0, N+1]$. However, even if $F_k, k=0, 1,\cdots, N$,
are nonsingular, the following stochastic difference equation
\begin{equation}
x_{k+1}=F_kx_k+G_kx_kw_k, \ x_{N+1}\in l^2_{{\cal F}_{N}}\label{eq
fdfdf}
\end{equation}
with terminal state given does not always admit an $l^2_{{\cal
F}_{k-1}}$-solution. For example, if we take $F_k=1$, $G_k=0$, and
the terminal state $x_2=w_1$ in (\ref{eq fdfdf}), then
$x_1=w_1\notin l^2_{{\cal F}_{0}}$, $x_0=w_1\notin l^2_{{\cal
F}_{-1}}$.
\end{example}

\begin{remark}\label{rem:2.2}\rm\
A class of backward stochastic difference equations arising from the study of
discrete stochastic maximum principle can be found in
\cite{linzhang}.
\end{remark}

To define and better understand the uniform detectability for system
(\ref{eq vcgtgh1}), we first give some lemmas.

\begin{lemma}\label{lem:2.2}\
(i) For system (\ref{eq vcgtgh1}),
$E\|x_l\|^2=E\|\phi_{l,k}x_k\|^2$ for $ l\ge k$, where it is assumed
that $\phi_{k,k}=I$, and $\phi_{l,k}$ is given by the following
iterative relation
\begin{equation}\label{eq 1}
\phi_{l,k}=\left[
\begin{array}{cc}
\phi_{l,k+1}F_k\\
\phi_{l,k+1}G_k
\end{array}
\right], \ l>k.
\end{equation}
(ii) $x_k\in l^2_{{\cal F}_{k-1}}$ if $F_i$ and $G_i$ are bounded
for $0\le i\le k-1$.
\end{lemma}

\begin{IEEEproof}\
(i)  can be shown by induction. For $k=l-1$, we have
\begin{align}
E\|x_{l}\|^2&=E[(F_{l-1}x_{l-1}+G_{l-1}x_{l-1}w_{l-1})^T(F_{l-1}x_{l-1}+G_{l-1}x_{l-1}w_{l-1})]\nonumber \\
&=E[x_{l-1}^T(F^T_{l-1}F_{l-1}+G_{l-1}^TG_{l-1})x_{l-1}]\nonumber\\
&=E\|\phi_{l,l-1}x_{l-1}\|^2. \nonumber
\end{align}
Hence, (\ref{eq 1}) holds for $k=l-1$. Assume that for $k=m<l-1$,
$E\|x_{l}\|^2=E\|\phi_{l,m}x_{m}\|^2$. Next, we prove
$E\|x_{l}\|^2=E\|\phi_{l,m-1}x_{m-1}\|^2.$ It can be seen that
\begin{align}
E\|x_{l}\|^2&=E[x_{m}^T\phi_{l,m}^T\phi_{l,m}x_{m}]\nonumber \\
&=E[(F_{m-1}x_{m-1}+G_{m-1}x_{m-1}w_{m-1})^T\phi_{l,m}^T\phi_{l,m}(F_{m-1}x_{m-1}+G_{m-1}x_{m-1}w_{m-1})]\nonumber\\
&=E[x_{m-1}^T(F_{m-1}^T\phi_{l,m}^T\phi_{l,m}F_{m-1}+G_{m-1}^T\phi_{l,m}^T\phi_{l,m}G_{m-1})x_{m-1}]\nonumber\\
&=E\|\phi_{l,m-1}x_{m-1}\|^2.\label{eq step1 in Lemma 1}\nonumber
\end{align}
This completes the proof of (i). And (ii) is obvious. The proof of this
lemma is complete.
\end{IEEEproof}

\begin{lemma}\label{lem:2.1}\
For system (\ref{eq vcgtgh1}), there holds
$\sum_{i=k}^lE\|y_i\|^2=E\|H_{l,k}x_k\|^2$ for $l\ge k\ge 0$, where
\begin{equation}
H_{l,k}=\left[
\begin{array}{cc}
H_k\\
(I_2\otimes H_{k+1})\phi_{k+1,k}\\(I_{2^2}\otimes
H_{k+2})\phi_{k+2,k}\\ \vdots \\ (I_{2^{l-k}}\otimes H_l)\phi_{l,k}
\end{array}
\right] \label{eq 2}
\end{equation}
with $H_{k,k}=H_k$ and $\phi_{j,k}(j=k+1,\cdots,l)$ given by
(\ref{eq 1}).
\end{lemma}

\begin{IEEEproof}\
We prove this lemma by induction. First, by a straight and simple computation,
the conclusion holds in the case of  $k=l,l-1$. Next, we assume that for
$k=m<l-1$, $\sum_{i=m}^{l}E\|y_i\|^2=E\|H_{l,m}x_m\|^2$ holds, then
it only needs to prove
$\sum_{i=m-1}^{l}E\|y_i\|^2=E\|H_{l,m-1}x_{m-1}\|^2$. It can be
verified that
\begin{align}\label{eq 9}
\sum_{i=m-1}^{l}E\|y_i\|^2&= \sum_{i=m}^{l}E\|y_i\|^2+E\|y_{m-1}\|^2\nonumber\\
&= E\|H_{l,m}x_m\|^2+E\|y_{m-1}\|^2\nonumber\\
&= E[x_m^TH_{l,m}^TH_{l,m}x_m]+E[x_{m-1}^TH_{m-1}^TH_{m-1}x_{m-1}]\nonumber\\
&= E[(F_{m-1}x_{m-1}+G_{m-1}x_{m-1}w_{m-1})^TH_{l,m}^TH_{l,m}(F_{m-1}x_{m-1}+G_{m-1}x_{m-1}w_{m-1})]\nonumber\\
&~~~+E[x_{m-1}^TH_{m-1}^TH_{m-1}x_{m-1}]\nonumber\\
&= E\left\{x_{m-1}^T\left[\begin{array}{c}
H_{m-1}\\
H_{l,m}F_{m-1}\\
H_{l,m}G_{m-1}
\end{array}
\right]^T
\left[\begin{array}{c}
H_{m-1}\\
H_{l,m}F_{m-1}\\
H_{l,m}G_{m-1}
\end{array}
\right]x_{m-1}\right\}.
\end{align}
By (\ref{eq 2}), it follows that
\begin{equation}\label{eq aa}
\left[\begin{array}{c}H_{m-1}\\H_{l,m}F_{m-1}\\
H_{l,m}G_{m-1}\end{array}\right]=\left[\begin{array}{c}H_{m-1}\\H_mF_{m-1}\\(I_2\otimes
H_{m+1})\phi_{m+1,m}F_{m-1}\\\vdots\\(I_{2^{t-m}}\otimes
H_l)\phi_{t,m}F_{m-1}\\H_mG_{m-1}\\(I_2\otimes
H_{m+1})\phi_{m+1,m}G_{m-1}\\\vdots\\(I_{2^{t-m}}\otimes
H_l)\phi_{l,m}G_{m-1}\end{array}\right].
\end{equation}
On the other hand, it can be deduced from (\ref{eq 1}) and (\ref{eq 2}) that
\begin{equation}\label{eq bb}
H_{l,m-1}=\left[\begin{array}{c}H_{m-1}\\(I_2\otimes
H_m)\left[\begin{array}{c}F_{m-1}\\G_{m-1}\end{array}\right]\\(I_{2^2}\otimes
H_{m+1})\left[\begin{array}{c}\phi_{m+1,m}F_{m-1}\\\phi_{m+1,m}G_{m-1}\end{array}\right]\\\vdots\\
(I_{2^{t-m+1}}\otimes
H_l)\left[\begin{array}{c}\phi_{l,m}F_{m-1}\\\phi_{l,m}G_{m-1}\end{array}\right]\end{array}\right].
\end{equation}
Combining (\ref{eq aa}) and (\ref{eq bb}) together results in
\begin{equation}\label{eq 12}
\left[\begin{array}{c}H_{m-1}\\H_{l,m}F_{m-1}\\
H_{l,m}G_{m-1}\end{array}\right]^T\left[\begin{array}{c}H_{m-1}\\H_{l,m}F_{m-1}\\
H_{l,m}G_{m-1}\end{array}\right]=H_{l,m-1}^TH_{l,m-1}.
\end{equation}
Hence, $\sum_{i=m-1}^{l}E\|y_i\|^2=E\|H_{l,m-1}x_{m-1}\|^2$. This
lemma is shown.

\end{IEEEproof}

Based on Lemmas~\ref{lem:2.2}--\ref{lem:2.1}, we are now in a position
to define the uniform detectability for system (\ref{eq vcgtgh1}).

\begin{definition}\label{def:2.2}\
System (\ref{eq vcgtgh1}) or $(F_k,G_k|H_k)$ is said to be uniformly
detectable if there exist integers $s, t\ge 0$, and positive constants
$d$, $b$ with $0\le d<1$ and $0<b<\infty$ such that whenever
\begin{equation}
E\|x_{k+t}\|^2=E\|\phi_{k+t,k}x_k\|^2\ge d^2 E\|x_k\|^2 , \label{eq
vhbjh1}
\end{equation}
there holds
\begin{equation}
 \sum_{i=k}^{k+s}
E\|y_i\|^2 =E\|H_{k+s,k}x_k\|^2\ge b^2 E\|x_k\|^2 , \label{eq
vhhjj1}
\end{equation}
where $k\in {\cal N}_0$, and $\phi_{k+t,k}$ and $H_{k+s,k}$ are the
same as defined in Lemma~\ref{lem:2.1}.
\end{definition}

Obviously, without loss of generality, in Definition~\ref{def:2.2}
we can assume that $t\le s$. By Lemmas~\ref{lem:2.2}--\ref{lem:2.1},
the uniform detectability of $(F_k,G_k|H_k)$ implies, roughly
speaking, that the state trajectory decays faster than the output energy
does. In what follows, ${\cal O}_{k+s,k}:=H_{k+s,k}^TH_{k+s,k}$ is
called an observability Gramian matrix, and $\phi_{l,k}$ a state
transition matrix from $x_k$ to $x_l$ of stochastic system (\ref{eq
vcgtgh1}). So (\ref{eq vhhjj1}) can be written as $E[x_k^T{\cal
O}_{k+s,k}x_k]\ge b^2 E\|x_k\|^2$. If $G_k\equiv 0$ for $k\ge 0$,
then system (\ref{eq vcgtgh1}) reduces to the following
deterministic system
\begin{equation}
\left\{
\begin{split}
x_{k+1}&=F_kx_k, \ x_0\in {\cal R}^n , \\
y_k&=H_kx_k ,
\end{split}
\right. \label{eq vdchcf}
\end{equation}
which was discussed in \cite{bdo, peters}.

Similarly, uniform observability can be defined as follows:

\begin{definition}\label{def:UO}\
System (\ref{eq vcgtgh1}) or $(F_k,G_k|H_k)$ is said to be uniformly
observable if there exist an integer $s\ge 0$ and a positive
constant $b>0$ such that
\begin{equation}
E\|H_{k+s,k}x_k\|^2\ge b^2 E\|x_k\|^2 \nonumber
\end{equation} holds for each initial condition $x_k\in l^2_{{\cal
F}_{k-1}}$, $k\in {\cal N}_0$.
\end{definition}

\begin{remark}\label{rem:UO}\rm\
Different from the uniform detectability concept, uniform
observability needs that any models (unstable and stable) should be
reflected by the output. This section concentrates on the uniform
detectability of system (\ref{eq vcgtgh1}), since it is weaker than
uniform observability. Uniform observability is also an important
concept, which will be further studied in the future.
\end{remark}

\begin{definition}\label{def:2.3}\
System (\ref{eq vcgtgh1}) is said to be exponentially stable in mean square
(ESMS) if there exist $\beta\ge 1$ and $\lambda\in (0,1)$ such that for any
$0\le k_0\le k<+\infty$, there holds
\begin{equation}
E\|x_k\|^2\le \beta E\|x_{k_0}\|^2 \lambda^{(k-k_0)}. \label{eq vcgv}
\end{equation}
\end{definition}

\begin{proposition}\label{prop:2.1}\
If system (\ref{eq vcgtgh1}) is ESMS, then for any bounded matrix sequence
$\{H_k\}_{k\ge 0}$, system (\ref{eq vcgtgh1}) is uniformly detectable.
\end{proposition}

\begin{IEEEproof}\
By Definition~\ref{def:2.3}, for any $k, t\ge 0$, we always have
\begin{equation}
E\|x_{k+t}\|^2=E\|\phi_{k+t,k}x_k\|^2\le \beta E\|x_{k}\|^2
\lambda^{t}, \ \beta>1, \ 0<\lambda<1. \label{eq vcgvff}
\end{equation}
By (\ref{eq vcgvff}), $\beta\lambda^t\to 0$ as $t\to\infty$. Set a
large $t_0>0$ such that $0\le d^2:=\beta\lambda^{t_0}<1$. Then, for
any fixed $t>t_0$, (\ref{eq vhbjh1}) holds only for $x_k=0$, which
makes (\ref{eq vhhjj1}) valid for any $s\ge t>t_0$ and $b>0$ with an
equality. So system (\ref{eq vcgtgh1}) is uniformly detectable.
\end{IEEEproof}

\begin{remark}\label{rem:2.3}\rm\
For system (\ref{eq vdchcf}), Definition~\ref{def:2.2} reduces to
Definition~2.1 in \cite{bdo}. It is easy to prove that uniform
detectability coincides with classical detectability of the linear
time-invariant system (\ref{eq vgthgh}).
\end{remark}

The following lemma will be used throughout this paper.

\begin{lemma}[see \cite{lguo}]\label{lem:2.3}\
For a nonnegative real sequence $\{s_k\}_{k\ge k_0}$, if there exist
constants $M_0\ge 1$, $\delta_0\in (0,1)$, and an integer $h_0>0$
such that $s_{k+1}\le M_0s_k$ and $\min_{k+1\le i\le k+h_0}s_i\le
\delta_0 s_k$, then
$$
s_k\le [M_0^{h_0}{\delta_0}^{-1}](\delta_0^{h_0})^{k-k_0}s_{k_0}, \
\forall k\ge k_0.
$$
\end{lemma}

The following proposition extends Lemma~2.2 in \cite{bdo}.

\begin{proposition}\label{prop:2.2}\
Suppose that $(F_k,G_k|H_k)$ is uniformly detectable, and $F_k$ and $G_k$
are uniformly bounded, i.e., $\|F_k\|\le M, \|G_k\|\le M, M>0$. Then
$\lim_{k\to\infty}E\|y_k\|^2=0$ implies $\lim_{k\to\infty}E\|x_k\|^2=0$.
\end{proposition}

\begin{IEEEproof}\
If there exists some integer $k_0$ such that for all $k\ge k_0$,
$E\|x_{k+t}\|^2=E\|\phi_{k+t,k}x_k\|^2<d^2E\|x_k\|^2$, then
$\min_{k+1\le i\le k+t}$ $E\|x_i\|^2<d^2E\|x_k\|^2$. Moreover,
$E\|x_{i+1}\|^2=E\|\phi_{i+1,i}x_i\|^2=E[x_i^T(F_i^TF_i+G_i^TG_i)x_i]\le
2M^2E\|x_i\|^2\le M_0E\|x_i\|^2$, where $M_0=\max\{2M^2,1\}\ge 1$.
By Lemma~\ref{lem:2.3}, not only does
$\lim_{k\to\infty}E\|x_k\|^2=0$, but also is system (\ref{eq
vcgtgh1}) ESMS. Otherwise, there exists a subsequence $\{k_i\}_{i\ge
0}$ such that $E\|\phi_{k_i+t,k_i}x_{k_i}\|^2\ge d^2E\|x_{k_i}\|^2$.
Now, for $k\in (k_i,k_{i+1})$, we write $k=k_i+1+t\alpha+\beta$ with
$\beta<t$, then
$$
E\|x_{k_i+1+\alpha t}\|^2\le d^\alpha E\|x_{k_i+1}\|^2,
$$
$$
E\|x_{k_i+1+\alpha t+\beta}\|^2\le (2M^2)^\beta E\|x_{k_i+1+\alpha
t}\|^2,
$$
$$
E\|x_{k_i+1}\|^2\le 2M^2 E\|x_{k_i}\|^2.
$$
Therefore, we have
\begin{align}
E\|x_k\|^2&= E\|x_{k_i+1+\alpha t+\beta}\|^2\le (2M^2)^\beta d^\alpha E\|x_{k_i+1}\|^2\nonumber\\
&\le (2M^2)^{\beta+1}d^\alpha E\|x_{k_i}\|^2. \label{eq vhvhbh}
\end{align}
Obviously, in order to show $\lim_{k\to\infty} E\|x_k\|^2=0$, we only need
to show $\lim_{k_i\to\infty} E\|x_{k_i}\|^2=0$. If it is not so, then there
are a subsequence $\{n_i\}_{i\ge 0}$ of $\{k_i\}_{i\ge 0}$ and $\varsigma>0$,
such that $E\|x_{n_i}\|^2>\varsigma$,
$E\|\phi_{n_i+t,n_i}x_{n_i}\|^2\ge d^2E\|x_{n_i}\|^2$. By Definition~\ref{def:2.2},
\begin{equation}
\sum_{i=n_i}^{n_i+s} E\|y_i\|^2=E[x_{n_i}^T{\cal
O}_{n_i+s,n_i}x_{n_i}]\ge b^2E\|x_{n_i}\|^2>b^2\varsigma. \label{eq vhhbb}
\end{equation}
Taking $n_i\to \infty$ in (\ref{eq vhhbb}), we have $0>b^2\varsigma>0$, which
is a contradiction. Hence, the proof is complete.
\end{IEEEproof}

In the remainder of this section, we will prove the output feedback
invariance for uniform detectability. Consider the following LDTV
stochastic control system
\begin{equation}
\left\{
\begin{split}
x_{k+1}&=(F_kx_k+M_ku_k)+(G_kx_k+N_ku_k)w_k, \\
y_k&=H_kx_k, \ k=0,1,2,\cdots .
\end{split}
\right. \label{eq vcgtgh1dd}
\end{equation}
Applying an output feedback control law $u_k=K_ky_k$ to (\ref{eq vcgtgh1dd})
yields the following closed-loop system
\begin{equation}
\left\{
\begin{split}
x_{k+1}&=(F_k+M_kK_kH_k)x_k+(G_k+N_kK_kH_k)x_kw_k, \\
y_k&=H_kx_k, \ k=0,1,2,\cdots .
\end{split}
\right. \label{eq vcgfdef}
\end{equation}

\begin{theorem}\label{thm:2.1}\
If $(F_k,G_k|H_k)$ is uniformly detectable, then so is
$(F_k+M_kK_kH_k, G_k+N_kK_kH_k|H_k)$.
\end{theorem}

\begin{IEEEproof}\
By Lemma~\ref{lem:2.1}, the observability Gramian for system (\ref{eq vcgfdef})
is $\bar {\cal O}_{k+s,k}={\bar H}^T_{k+s,k}{\bar H}_{k+s,k}$, where
$$
{\bar H}_{k+s,k}=\left[
\begin{array}{cc}
H_k\\
(I_2\otimes H_{k+1}){\bar \phi}_{k+1,k}\\(I_{2^2}\otimes
H_{k+2}){\bar \phi}_{k+2,k}\\ \vdots \\ (I_{2^{s}}\otimes
H_{k+s}){\bar \phi}_{{k+s},k}
\end{array}
\right], \ \bar \phi_{k+i,k}=\left[
\begin{array}{cc}
\bar \phi_{k+i,k+1}\bar F_k\\
\bar \phi_{k+i,k+1}\bar G_k \end{array} \right], \ i=1,\cdots, s.
$$
$$
{\bar F}_j=F_j+M_jK_jH_j, \ {\bar G}_j=G_j+N_jK_jH_j, \ j=k, k+1,
\cdots, k+s.
$$
To prove that $({\bar F}_k, {\bar G}_k|H_k)$ is uniformly
detectable, it suffices to show that there are constants ${\bar
b}>0$, $0<\bar d<1$, $s, t\ge 0$ such that for $\xi\in l^2_{{\cal
F}_{k-1}}$, $k\in {\cal N}_0$, whenever
\begin{equation}
E[x_k^T\bar {\cal O}_{k+s,k}x_k]<{\bar b}^2 E\|x_k\|^2 , \label{eq
vhhjj}
\end{equation}
we have
\begin{equation}
E\|{\bar \phi}_{k+t,k}x_k\|^2<{\bar d}^2 E\|x_k\|^2 . \label{eq
vhbjh}
\end{equation}

It is easy to show
$$
{\bar H}_{k+s,k}=Q_{k+s,k}H_{k+s,k}, \ Q_{k+s,k}=\left[
\begin{array}{cccccc}
I & 0 & \cdots & 0\\
{*} & I & \cdots & 0\\
\vdots & \vdots & {}& \vdots\\
{*} & * &\cdots & I
\end{array}
\right].
$$
where * represents terms involving $H_i$, $M_i$, $K_i$ and $N_i$,
$i=k, k+1, \cdots, k+s$. Hence, for any $x_k\in l^2_{{\cal
F}_{k-1}}$,
\begin{eqnarray}
\rho E[x_k^T {\cal O}_{k+s,k}x_k] \le E[x_k^T\bar {\cal
O}_{k+s,k}x_k] \le \varrho E[x_k^T {\cal O}_{k+s,k}x_k] , \label{eq
cvghv}
\end{eqnarray}
where $\rho=\lambda_{\min} (Q_{k+s,k}^TQ_{k+s,k})>0$,
$\varrho=\lambda_{\max} (Q_{k+s,k}^TQ_{k+s,k})>0$. In addition, by observation,
for any $l>k\ge 0$,
$$
{\bar \phi}_{l,k}=\phi_{l,k}+R_{l,k} H_{l,k} ,
$$
where $R_{l,k}$ is a matrix involving $H_i$, $M_i$, $K_i$ and $N_i$,
$i=k, k+1, \cdots, l-1$. If we take $0<\bar b\le {\sqrt \rho} b$,
then it follows from (\ref{eq cvghv}) that $E[x_k^T {\cal
O}_{k+s,k}x_k]<\frac 1 \rho E[x_k^T\bar {\cal O}_{k+s,k}x_k]\le
\frac {\bar b^2}{\rho} E\|x_k\|^2\le b^2E\|x_k\|^2$. By the uniform
observability of $(F_k,G_k|H_k)$, it follows that
\begin{align}
E\|{\bar \phi}_{k+t,k}x_k\|^2
&= E\|\phi_{k+t,k}x_k+R_{k+t,k}H_{k+t,k}x_k\|^2\nonumber\\
&\le 2E\|{\phi}_{k+t,k}x_k\|^2+2\mu^2 E\|H_{k+t,k}x_k\|^2\nonumber\\
&\le 2d^2E\|x_k\|^2+2\mu^2E[x_k^T {{\cal O}}_{k+s,k}x_k]\nonumber\\
&\le \big(2d^2+2\mu^2\frac {\bar b^2} {\rho}\big)E\|x_k\|^2\nonumber\\
&= \bar d E\|x_k\|^2 ,
\end{align}
where $\mu=\sup_{k} \|R_{k+t,k}\|$, $\bar d=2d^2+2\mu^2\frac {\bar
b^2} {\rho}$. If we take $\bar b$ to be sufficiently small, then
$\bar d<1$, which yields the uniform detectability of $({\bar F}_k,
{\bar G}_k|H_k)$. Hence, the proof of this theorem is complete.
\end{IEEEproof}

Theorem~\ref{thm:2.1} reveals that the output feedback does not change
uniform detectability.

\begin{example}\label{exm:2.4}\rm\
For simplicity, we set $s=1$. Then it can be computed that
$$
{\bar H}_{k+1,k}=\left[
\begin{array}{cc}
H_k\\
(I_2\otimes H_{k+1}){\bar\phi}_{k+1,k}
\end{array}
\right]=\left[
\begin{array}{cc}
H_k \\
{H_{k+1}}(F_k+M_kK_kH_k)\\
{H_{k+1}}(G_k+N_kK_kH_k)
\end{array}
\right],
$$
$$
{H}_{k+1,k}=\left[
\begin{array}{cc}
H_k\\
(I_2\otimes H_{k+1}){\phi}_{k+1,k}
\end{array}
\right]=\left[
\begin{array}{cc}
H_k \\
{H_{k+1}}F_k\\
{H_{k+1}}G_k
\end{array}
\right].
$$
Obviously,
$$
Q_{k+1,k}=\left[
\begin{array}{ccc}
I & 0 & 0 \\
{H_{k+1}}M_kK_k & I & 0\\
{H_{k+1}}N_kK_k & 0 & I
\end{array}
\right].
$$
\end{example}

\begin{example}\label{exm:2.5}\rm\
By definition, we have
$$
{\bar\phi}_{k+1,k}=\left[
\begin{array}{ccc}
F_k+M_kK_kH_k\\
G_k+N_kK_kH_k
\end{array}
\right], \ \ {\phi}_{k+1,k}=\left[
\begin{array}{ccc}
F_k\\
G_k
\end{array}
\right].
$$
Hence, ${\bar\phi}_{k+1,k}={\phi}_{k+1,k}+R_{k+1,k}H_{k+1,k}$ with
$R_{k+1,k}=\left[
\begin{array}{ccc}
M_kK_k & 0 & 0\\
N_kK_k & 0 & 0
\end{array}
\right].
$
\end{example}

\subsection{Lyapunov-Type Theorems under Uniform Detectability\label{sec:TLTT-1}}

In the following, we will further study the following time-varying
GLE
\begin{equation}
-P_k+F_k^TP_{k+1}F_k+G_k^TP_{k+1}G_k+H_k^TH_k=0, \ k=0,1,2,\cdots,
\label{eq vhvcxf}
\end{equation}
under uniform detectability. The aim of this subsection is to extend
the classical Lyapunov theorem to GLE (\ref{eq vhvcxf}). To study
(\ref{eq vhvcxf}), we first introduce the following finite time
backward difference equation
\begin{equation}
\left\{
\begin{array}{l}
-P_{k,T}+F_{k}^TP_{k+1,T}F_k+G_k^TP_{k+1,T}G_k+H_k^TH_k=0,\\
P_{T,T}=0, \ k=0, 1, \cdots, T-1; \ T\in {\cal N}_1.
\end{array}
\right. \label{eq vceqwe}
\end{equation}
Obviously, equation (\ref{eq vceqwe}) has nonnegative definite
solutions $P_{k,T}\ge 0$.

\begin{propsub}\label{prop:4.1.1}\
$P_{k,T}$ is monotonically increasing with respect to $T$, i.e., for
any $k_0\le T_1\le T_2<+\infty$,
$$
P_{k_0,T_1}\le P_{k_0,T_2}, \ k_0\in \{0,1, \cdots, T_1\}.
$$
\end{propsub}

\begin{IEEEproof}\
Obviously, $P_{k,T_1}$ and $P_{k,T_2}$ solve
\begin{equation}
\left\{
\begin{array}{l}
-P_{k,T_1}+F_{k}^TP_{k+1,T_1}F_k+G_k^TP_{k+1,T_1}G_k+H_k^TH_k=0,\\
P_{T_1,T_1}=0, \ k=0, 1, \cdots, T_1-1 ,
\end{array}
\right. \label{eq vdeffe}
\end{equation}
and
\begin{equation}
\left\{
\begin{array}{l}
-P_{k,T_2}+F_{k}^TP_{k+1,T_2}F_k+G_k^TP_{k+1,T_2}G_k+H_k^TH_k=0,\\
P_{T_2,T_2}=0, \ k=0, 1, \cdots, T_2-1 ,
\end{array}
\right. \label{eq vdeffess}
\end{equation}
respectively. Consider the following LDTV stochastic system with a
deterministic initial state $x_{k_0}$:
\begin{equation}
\left\{
\begin{array}{l}
x_{k+1}=F_kx_k+G_kx_kw_k, \\
x_{k_0}\in {\cal R}^n , \ k=k_0, k_0+1, \cdots.
\end{array}
\right. \label{eq vcscdc}
\end{equation}
Associated with (\ref{eq vcscdc}), in view of (\ref{eq vdeffe}), we
have
\begin{align}
\sum_{k=k_0}^{T_1-1} E[x_k^TH_k^TH_kx_k]
&= \sum_{k=k_0}^{T_1-1} E[x_k^TH_k^TH_kx_k+x_{k+1}^TP_{k+1,T_1}x_{k+1}-x^T_kP_{k,T_1}x_k]\nonumber\\
&~~~ +x_{k_0}^TP_{k_0,T_1}x_{k_0}-E[x^T_{T_1}P_{T_1,T_1}x_{T_1}]\nonumber\\
&= \sum_{k=k_0}^{T_1-1} E[x_k^T(-P_{k,T_1}+F_{k}^TP_{k+1,T_1}F_k+G_k^TP_{k+1,T_1}G_k+H_k^TH_k)x_k]\nonumber\\
&~~~ +x_{k_0}^TP_{k_0,T_1}x_{k_0}\nonumber\\
&= x_{k_0}^TP_{k_0,T_1}x_{k_0}. \label{eq hvbhcdgt}
\end{align}
Similarly,
\begin{eqnarray}
\sum_{k=k_0}^{T_2-1}
E[x_k^TH_k^TH_kx_k]=x_{k_0}^TP_{k_0,T_2}x_{k_0}. \label{eq
hvbhcdgt1}
\end{eqnarray}
From (\ref{eq hvbhcdgt})-(\ref{eq hvbhcdgt1}), it follows that
\begin{eqnarray}
0 \le\sum_{k=k_0}^{T_1-1}
E[x_k^TH_k^TH_kx_k]=x_{k_0}^TP_{k_0,T_1}x_{k_0} \le
\sum_{k=k_0}^{T_2-1}
E[x_k^TH_k^TH_kx_k]=x_{k_0}^TP_{k_0,T_2}x_{k_0}. \label{eq vbhbjhb}
\end{eqnarray}
The above expression holds for any $x_{k_0}\in {\cal R}^n$, which
yields $P_{k_0,T_1}\le P_{k_0,T_2}$. Thus, the proof is complete.
\end{IEEEproof}

\begin{propsub}\label{prop:4.1.2}\
If system (\ref{eq vcgtgh1}) is ESMS, and $H_k$ is uniformly bounded
(i.e., there exists $M>0$ such that $\|H_k\|\le M$, $\forall k\in
{\cal N}_0$), then the solution $P_{k,T}$ of (\ref{eq vceqwe}) is
uniformly bounded for any $T\in {\cal N}_1$ and $k\in [0,T]$.
\end{propsub}

\begin{IEEEproof}\
By (\ref{eq hvbhcdgt}), for any deterministic $x_k\in {\cal R}^n$,
we have
\begin{align*}
x_{k}^TP_{k,T}x_{k}
&= \sum_{i=k}^{T-1} E[x_i^TH_i^TH_ix_i]\le \sum_{i=k}^{\infty} E[x_i^TH_i^TH_ix_i]\\
&\le M^2\|x_k\|^2\beta
\sum_{i=k}^{\infty}\lambda^{(i-k)}=M^2\|x_k\|^2 \beta
\frac{1}{1-\lambda} ,
\end{align*}
which leads to that $0\le P_{k,T}\le \frac {\beta M^2}{1-\lambda} I$
since $x_{k}$ is arbitrary. Hence, the proof is complete.
\end{IEEEproof}

Combining Proposition~\ref{prop:4.1.1} with
Proposition~\ref{prop:4.1.2} yields that $P_k:=\lim_{T\to\infty}
P_{k,T}$ exists, which is a solution of (\ref{eq vhvcxf}). Hence, we
obtain the following Lyapunov-type theorem.

\begin{thmsub}[Lyapunov-Type Theorem]\label{thm:4.1.1}\
If system (\ref{eq vcgtgh1}) is ESMS and $\{H_k\}_{k\in {\cal N}_0}$
is uniformly bounded, then (\ref{eq vhvcxf}) admits a unique
nonnegative definite solution $\{P_k\}_{k\in {\cal N}_0}$.
\end{thmsub}

The converse of Theorem~\ref{thm:4.1.1} still holds.

\begin{thmsub}[Lyapunov-Type Theorem]\label{thm:4.1.2}\
Suppose that $(F_k,G_k|H_k)$ is uniformly detectable and $F_k$ and
$G_k$ are uniformly bounded with an upper bound $M>0$. If there is a
bounded nonnegative definite symmetric matrix sequence
$\{P_k\}_{k\ge 0}$ solving GLE (\ref{eq vhvcxf}), then system
(\ref{eq vcgtgh1}) is ESMS.
\end{thmsub}

\begin{IEEEproof}\
For system (\ref{eq vcgtgh1}), we take a Lyapunov function as
$$
V_k(x)=x^T(P_k+\varepsilon I)x ,
$$
where $\varepsilon>0$ is to be determined. For simplicity, in the
sequel, we let $V_k:=V_k(x_k)$. It is easy to compute
\begin{align}
EV_k-EV_{k+1}
&= E[x_k^T(P_k+\varepsilon I)x_k]-E[x_{k+1}^T(P_{k+1}+\varepsilon I)x_{k+1}]\nonumber\\
&= E[x_k^T(P_k+\varepsilon I)x_k]-E[(F_kx_k+G_kx_kw_k)^T(P_{k+1}
+\varepsilon I)
(F_kx_k+G_kx_kw_k)]\nonumber\\
&= E[x_k^T(P_k-F^T_kP_{k+1}F_k-G_k^TP_{k+1}G_k)x_k]+\varepsilon
E[x_k^T(I-F_k^TF_k-G_k^TG_k)x_k]\nonumber\\
&= E\|y_k\|^2+\varepsilon E[x_k^T(I-F_k^TF_k-G_k^TG_k)x_k]\nonumber\\
&= E\|y_k\|^2+\varepsilon E\|x_k\|^2-\varepsilon E\|x_{k+1}\|^2.
\label{eq vhbh}
\end{align}
Identity (\ref{eq vhbh}) yields
\begin{align}
EV_k-EV_{k+s+1}
&= [EV_k-EV_{k+1}]+[EV_{k+1}-EV_{k+2}]+\cdots+[EV_{k+s}-EV_{k+s+1}]\nonumber\\
&= \sum_{i=k}^{k+s}E\|y_i\|^2+\varepsilon E\|x_k\|^2-\varepsilon
E\|x_{k+s+1}\|^2. \label{eq fsdfd}
\end{align}
When $\sum_{i=k}^{k+s}E\|y_i\|^2\ge b^2E\|x_k\|^2$, we first note
that
\begin{align}
E\|x_{k+s+1}\|^2
&= E\{x^T_{k+s}(F_{k+s}^TF_{k+s}+G_{k+s}^TG_{k+s})x_{k+s}\}\nonumber\\
&\le 2M^2E\|x_{k+s}\|^2\le (2M^2)^2 E\|x_{k+s-1}\|^2\le \cdots \nonumber\\
&\le (2M^2)^{s+1}E\|x_k\|^2. \label{eq hdbdj}
\end{align}
Then, by (\ref{eq fsdfd}), we still have
\begin{align}
EV_k-EV_{k+s+1}
&\ge b^2E\|x_k\|^2+\varepsilon E\|x_k\|^2-\varepsilon (2M^2)^{s+1}E\|x_k\|^2\nonumber\\
&= [b^2+\varepsilon-\varepsilon (2M^2)^{s+1}]E\|x_k\|^2. \label{eq
vhvhn}
\end{align}
From (\ref{eq vhvhn}), it readily follows that
\begin{align}
EV_{k+s+1}
&\le EV_k-\{b^2+\varepsilon [1-(2M^2)^{s+1}]\}E\|x_k\|^2\nonumber\\
&\le
\left\{1-\frac{[b^2+\varepsilon[1-(2M^2)^{s+1}]]}{\lambda_{max}(P_k+\varepsilon
I)}\right\}EV_k. \label{eq cvgvg}
\end{align}
Considering that $\{P_k\ge 0\}_{k\in {\cal N}_0}$ is uniformly
bounded, if $\varepsilon$ is taken to be sufficiently small, then
there must exist a $\delta\in (0,1)$ such that
\begin{eqnarray}
EV_{k+s+1}\le \delta EV_k. \label{eq vvhbb}
\end{eqnarray}
When $\sum_{i=k}^{k+s}E\|y_i\|^2\le b^2E\|x_k\|^2$, by uniform
detectability we have $E\|x_{k+t}\|^2\le d^2 E\|x_k\|^2$. From
(\ref{eq fsdfd}), it follows that
\begin{eqnarray}
EV_k-EV_{k+t} \ge \varepsilon E\|x_k\|^2-\varepsilon d^2 E\|x_k\|^2
=\varepsilon (1-d^2)E\|x_k\|^2. \label{eq dshbd}
\end{eqnarray}
Similarly, we can show that there exists a constant $\delta_1\in
(0,1)$ such that
\begin{eqnarray}
EV_{k+t}\le \delta_1 EV_k. \label{eq vcdcdv}
\end{eqnarray}
Set $\delta_0:=\max\{\delta,\delta_1\}$, in view of (\ref{eq vvhbb})
and (\ref{eq vcdcdv}), we have
\begin{eqnarray}
\min_{k+1\le i\le k+s+1}{EV_i}\le \delta_0 EV_k,\ \forall k\ge 0.
\label{eq vhvh}
\end{eqnarray}
From identity (\ref{eq vhbh}), we know
\begin{eqnarray}
EV_{k+1} \le EV_k+\varepsilon E\|x_{k+1}\|^2\le EV_k+\varepsilon
EV_{k+1}. \label{eq vhvhhk}
\end{eqnarray}
Taking $0<\varepsilon<1$ in (\ref{eq vhvhhk}), it is easy to derive
that there exists a positive constant $M_0\ge 1$ satisfying
\begin{eqnarray}
EV_{k+1}\le M_0 EV_k, \ \forall k\ge 0. \label{eq vhvhddd}
\end{eqnarray}
Applying Lemma~\ref{lem:2.3} with $s_k=EV_k$, $h_0=s+1$,
$\beta=[M_0^{h_0}{\delta_0}^{-1}]$, $\lambda=\delta_0^{h_0}$, it
follows that
$$
EV_{k}\le \beta \lambda^{(k-k_0)}EV_{k_0} \le
\lambda_{\max}(P_k+\varepsilon I) \beta
\lambda^{(k-k_0)}E\|x_{k_0}\|^2 ,
$$
which implies that system (\ref{eq vcgtgh1}) is ESMS due to the fact
that $\{P_k\}_{k\ge 0}$ is uniformly bounded.
\end{IEEEproof}

The above theorem directly yields the following result.

\begin{corsub}\label{coro:4.1.3}\
Suppose that there exists $\epsilon>0$ such that $H_k^TH_k>\epsilon
I$ for $k\in {\cal N}_{0}$. Additionally, if there is a uniformly
bounded symmetric matrix sequence $\{P_k\ge 0\}_{k\ge 0}$ solving
GLE (\ref{eq vhvcxf}), then system (\ref{eq vcgtgh1}) is ESMS.
\end{corsub}

%%%%%%%%%%%%%%%%%%%%%%%%%%%%%%%%%%%%%%%%%%%%%%%%%%%%%%%%%%%%%%%%%%%%%%%%%%%%%%%%
\section{Exact Detectability and Related Lyapunov-type Theorems \label{sec:EDEO}}
%%%%%%%%%%%%%%%%%%%%%%%%%%%%%%%%%%%%%%%%%%%%%%%%%%%%%%%%%%%%%%%%%%%%%%%%%%%%%%%%

We recall that for the linear time-invariant system
\begin{equation}
\left\{
\begin{split}
x_{k+1}&=Fx_k+Gx_kw_k, \ x_0\in {\cal R}^n\\
y_k&=Hx_k, \ k=0,1,2,\cdots
\end{split}
\right. \label{eq vcgdsdsf}
\end{equation}
its exact observability was defined in \cite{zhang1,lizao}, while
the same definition for linear continuous-time time-invariant It\^o
systems was given in \cite{c3dd}. For the LDTV stochastic system
(\ref{eq vcgtgh1}), the complete observability that is different
from the uniform observability \cite{c2ddss} was defined in
\cite{xsdsffdf}. In this section, we will study exact detectability
of the stochastic system (\ref{eq vcgtgh1}), from which it can be found
that there are some essential differences between the time-varying
and time-invariant systems. In addition, Lyapunov-type theorems are
also presented.

%%%%%%%%%%%%%%%%%%%%%%%%%%%%%%%%%%%%%%%%%%%%%%%%%%%%%%%%%%%%%%%%%%%%%%%%%%%%%%%%
\subsection{Exact Detectability\label{sec:EDEO-1}}

We first give several definitions.

\begin{defsub}\label{def:3.1.1}\
For system (\ref{eq vcgtgh1}), $x_{k_0}\in {l}^2_{{\cal F}_{k_0-1}}$
is called a $k_0^\infty$-unobservable state if $y_k\equiv 0$ a.s.
for $k\in [{k_0},\infty)$, and $x_{k_0}\in {l}^2_{{\cal F}_{k_0-1}}$
is called a $k_0^{s_0}$-unobservable state if $y_k\equiv 0$ a.s. for
$k\in [k_0, k_0+s_0]$.
\end{defsub}

\begin{remsub}\label{rem:3.1.1}\rm\
From Definition~\ref{def:3.1.1}, we point out the following obvious facts:
(i) If $x_{k_0}$ is a $k_0^\infty$-unobservable state, then for any $s_0\ge 0$,
it must be a $k_0^{s_0}$-unobservable state; (ii) If $x_{k_0}$ is a
$k_0^{s_1}$-unobservable state, then for any $0\le s_0\le s_1$, it must be
a $k_0^{s_0}$-unobservable state.
\end{remsub}

\begin{exmsub}\label{exm:3.1.1}\rm\
In system (\ref{eq vcgtgh1}), if we take $H_k\equiv 0$ for $k\ge
k_0$, then any state $x_{k_0}\in {l}^2_{{\cal F}_{k_0-1}}$ is a
$k_0^\infty$-unobservable state. For any $k_0\ge 0$, $x_{k_0}=0$ is
a trivial $k_0^\infty$-unobservable state.
\end{exmsub}

Different from the linear time-invariant system (\ref{eq vcgdsdsf}),
even if $x_{k_0}=\zeta$ is a $k_0^{\infty}$-unobservable state,
$x_{k_1}=\zeta$ may not be a $k_1^{s_1}$-unobservable state for any
$s_1\ge 0$, which is seen from the next example.

\begin{exmsub}\label{exm:3.1.2}\rm\
Consider the deterministic linear time-varying system with $G_k=0$ and
$$
H_k=F_k=\left\{
\begin{array}{l}
\left[
\begin{array}{ccc}
1 & 0 \\
0 & 0
\end{array}
\right], \ \mbox{if} \ k \ \mbox{is} \ \mbox{even},\\
\left[
\begin{array}{ccc}
0 & 0\\
0 & 1
\end{array}
\right], \
 \mbox{if} \ k \ \mbox{is} \ \mbox{odd}.
 \end{array}
 \right.
 $$
Obviously, $x_0=\left[
\begin{array}{ccc}
0 \\
1
\end{array}
\right]$ is a $0^\infty$-unobservable state due to $y_k=0$ for $k\ge 0$,
but $x_1=\left[
\begin{array}{ccc}
0 \\
1
\end{array}
\right]$ is not a $1^{s_1}$-unobservable state for any $s_1\ge 0$
due to $y_1=H_1\left[
\begin{array}{ccc}
 0 \\
 1
\end{array}
\right]\ne 0$, let alone $1^{\infty}$-unobservable state.
\end{exmsub}

\begin{defsub}\label{def:3.1.2}\
System (\ref{eq vcgtgh1}) is called $k_0^\infty$-exactly detectable if all
$k_0^\infty$-unobservable state $\xi$ is exponentially stable, i.e., there
are constants $\beta\ge 1$, $0<\lambda<1$ such that
\begin{equation}
E\|x_k\|^2\le \beta E\|\xi\|^2 \lambda^{(k-k_0)}, \ \
 \forall k\ge k_0. \label{eq vcgvass}
\end{equation}
Similarly, system (\ref{eq vcgtgh1}) is called $k_0^{s_0}$-exactly detectable
if (\ref{eq vcgvass}) holds for all $k_0^{s_0}$-unobservable state $\xi$.
\end{defsub}

\begin{defsub}\label{def:3.1.3}\
System (\ref{eq vcgtgh1}) (or $(F_k,G_k|H_k)$) is said to be
${\cal K}^\infty$-exactly detectable if it is $k^\infty$-exactly detectable
for any $k\ge 0$. If there exists a nonnegative integer sequence
$\{s_k\}_{k\ge 0}$ with the upper limit
$\overline {\lim}_{k\to\infty}{s_k}=+\infty$ such that system
(\ref{eq vcgtgh1}) is $k^{s_k}$-exactly detectable, i.e., for
any $k^{s_k}$-unobservable state $\xi_k$,
$$
E\|x_t\|^2\le \beta E\|\xi_{k}\|^2 \lambda^{(t-k)}, \ \beta\ge 1, \
0<\lambda<1,\ t\ge k ,
$$
then system (\ref{eq vcgtgh1}) (or $(F_k,G_k|H_k)$) is said to be weakly
finite time or ${\cal K}^{WFT}$-exactly detectable. If \,
$\overline{\lim}_{k\to\infty}{s_k}<+\infty$, then system (\ref{eq vcgtgh1})
(or $(F_k,G_k|H_k)$) is said to be finite time or ${\cal K}^{FT}$-exactly
detectable.
\end{defsub}

A special case of ${\cal K}^{FT}$-exact detectability is the so-called
${\cal K}^{N}$-exact detectability, which will be used to study GLEs.

\begin{defsub}\label{def:3.1.4}\
If there exists an integer $N\ge 0$ such that for any time $k_0\in
[0,\infty)$, system (\ref{eq vcgtgh1}) (or $(F_k,G_k|H_k)$) is
$k_0^{N}$-exactly detectable, then system (\ref{eq vcgtgh1}) (or
$(F_k,G_k|H_k)$) is said to be ${\cal K}^{N}$-exactly detectable.
\end{defsub}

From Definitions~\ref{def:3.1.3}--\ref{def:3.1.4}, we have the following
inclusion relation
\begin{align*}
& {\cal K}^{N}\mbox {-exact detectability} \ \Longrightarrow \
{\cal K}^{FT}\mbox{-exact detectability} \\
\Longrightarrow \ & {\cal K}^{WFT}\mbox{-exact detectability} \
\Longrightarrow \ {\cal K}^{\infty}\mbox {-exact detectability}.
\end{align*}

In this paper, we will mainly use ${\cal K}^{\infty}$- and
${\cal K}^{N}$-exact detectability. Obviously, ${\cal K}^{N}$-exact
detectability implies ${\cal K}^{\infty}$-exact detectability, but the
converse is not true. We present the following examples to illustrate various
relations among several definitions on detectability. For illustration
simplicity, we only take the concerned examples to be deterministic.

\begin{exmsub}\label{exm:3.1.3}\rm\
In system (\ref{eq vdchcf}), we take $F_k=1$ for $k\ge 0$, and
$$
H_{k}=\left\{
\begin{array}{l}
1, \ \mbox{for} \ k=n^2, \ n=1, 2,\cdots,\\
0, \ \mbox{otherwise}.
 \end{array}
 \right.
$$
In this case, system (\ref{eq vdchcf}) (or $(F_k|H_k)$) is ${\cal
K}^{\infty}$-exactly detectable, and the zero vector is the unique
$k^{\infty}$-unobservable state. $(F_k|H_k)$ is also ${\cal
K}^{WFT}$-exactly detectable, where $s_k=k^2-k\to\infty$. However,
$(F_k|H_k)$ is not ${\cal K}^{FT}$-exactly detectable, and,
accordingly, is not ${\cal K}^{N}$-exactly detectable for any $N\ge
0$.
\end{exmsub}

\begin{exmsub}\label{exm:3.1.4}\rm\
In system (\ref{eq vdchcf}), if we take $F_k=1$ and $H_k=\frac 1 {k}$ for
$k\ge 0$, then $(F_k|H_k)$ is ${\cal K}^N$-exactly detectable for any
$N\ge 0$, but $(F_k|H_k)$ is not uniformly detectable. This is because for
any $t\ge 0$, $0\le d <1$ and $\xi\in {\cal R}$, we always have
$|\phi_{k+t,k}\xi|^2=|\xi|^2\ge d^2 |\xi|^2$. But there do not exist $b>0$
and $s\ge 0$ satisfying (\ref{eq vhhjj1}), because $\xi^T {\cal
O}_{k+s,k}\xi=|\xi|^2\sum_{i=k}^{k+s} \frac 1 {i^2}$ while
$\lim_{k\to\infty} \sum_{i=k}^{k+s} \frac 1 {i^2}=0$.
\end{exmsub}

\begin{exmsub}\label{exm:3.1.5}\rm\
In system (\ref{eq vdchcf}), if we take $F_k=1$ for $k\ge 0$, and
$H_ {2n}=1$ and $H_ {2n+1}=0$ for $n=0, 1, 2, \cdots$, then $(F_k|H_k)$ is
uniformly detectable and ${\cal K}^1$-exactly detectable, but it is not
${\cal K}^0$-exactly detectable.
\end{exmsub}

The following lemma is obvious.

\begin{lemsub}\label{lem:3.1.1}\
At any time $k_0$, $x_{k_0}=0$ is not only a $k_0^\infty$- but also a
$k_0^{s_0}$-unobservable state for any $s_0\ge 0$.
\end{lemsub}

By Lemma~\ref{lem:3.1.1}, if we let ${\Theta}^\infty_{k_0}$ denote the set of
all the $k_0^{\infty}$-unobservable states of system (\ref{eq vcgtgh1}) at
time $k_0$, then ${\Theta}^\infty_{k_0}$ is not empty. Furthermore, it is easy
to show that ${\Theta}^\infty_{k_0}$ is a linear vector space.

\begin{lemsub}\label{lem:3.1.2}\
For $k_0\in {\mathcal N}_0$, if there does not exist a nonzero
$\zeta\in {\mathcal R}^n$ such that $H_{k_0}\zeta=0$,
$(I_{2^{l-k_0}}\otimes H_{l})\phi_{l,k_0}\zeta=0$, $l=k_0+1, k_0+2,
\cdots$, then $y_k\equiv 0$ a.s. with $k\ge k_0$ implies $x_{k_0}=0$
a.s..
\end{lemsub}

\begin{IEEEproof}\
From $y_{k_0}\equiv 0$ a.s., it follows that
\begin{equation}
E[x_{k_0}^TH^T_{k_0}H_{k_0}x_{k_0}]=0. \label{eq gfff}
\end{equation}
From $y_l\equiv 0$ a.s., $l=k_0+1$, $\cdots$, it follows from
Lemma~\ref{lem:2.1} that
\begin{equation}
E\big[x_{k_0}^T\phi_{l,k_0}^T( I_{2^{l-k_0}} \otimes H^T_{l})(
I_{2^{l-k_0}} \otimes H_{l})\phi_{l,k_0}x_{k_0}\big]=0. \label{eq
vghvh}
\end{equation}
Let $R_{k_0}=E[x_{k_0}x_{k_0}^T]$, rank$R_{k_0}=r$. When $r=0$, this
implies $x_{k_0}=0$ a.s., and this lemma is shown. For $1\le r\le
n$, by the result of \cite{rantz}, there are real nonzero vectors
$z_1, z_2, \cdots, z_r$ such that $R_{k_0}=\sum_{i=1}^r z_iz_i^T$.
By (\ref{eq gfff}), we have
\begin{align}
E[x_{k_0}^TH^T_{k_0}H_{k_0}x_{k_0}]
&= \mathrm{trace}E[H^T_{k_0}H_{k_0}x_{k_0}x_{k_0}^T]\nonumber\\
&= \mathrm{trace}\{H^T_{k_0}H_{k_0}E[x_{k_0}x_{k_0}^T]\}\nonumber\\
&= \mathrm{trace}\{H^T_{k_0}H_{k_0} \sum_{i=1}^r z_iz_i^T\} \nonumber\\
&= \sum_{i=1}^r [z^T_iH^T_{k_0}H_{k_0} z_i]=0 , \label{eq bjbjb}
\end{align}
which gives $H_{k_0} z_i=0$ for $i=1,2, \cdots, r$. Similarly,
(\ref{eq vghvh}) yields
$$
( I_{2^{l-k_0}} \otimes H_{l})\phi_{l,k_0}z_i=0, \ i=1, 2, \cdots,
r.
$$
According to the given assumptions, we must have $z_i=0$, $i=1, 2,
\cdots, r$, which again implies $x_{k_0}=0$ a.s..
\end{IEEEproof}

By Lemma~\ref{lem:3.1.2}, it is known that under the conditions of
Lemma~\ref{lem:3.1.2}, $x_{k_0}=0$ is the unique $k_0^{\infty}$-unobservable
state, i.e., ${\Theta}^\infty_{k_0}=\{0\}$.

\begin{lemsub}\label{lem:3.1.3}\
Uniform detectability implies ${\cal K}^{\infty}$-exact detectability.
\end{lemsub}

\begin{IEEEproof}\
For any ${k_0}^{\infty}$-unobservable state $x_{k_0}=\xi$, by
Definition~\ref{def:2.2} and Definition~\ref{def:3.1.3}, we must have
$E\|\phi_{k+t,k}x_k\|^2<d^2E\|x_k\|^2$ or $x_k\equiv 0$ for $k\ge k_0$;
otherwise, it will lead to a contradiction since
$$
0=\sum_{i=k}^{k+s} E\|y_i\|^2\ge b\|x_k\|^2>0.
$$
Under any case, the following system
\begin{equation}
\left\{
\begin{split}
x_{k+1}&=F_kx_k+G_kx_kw_k, \\
x_{k_0}&=\xi \in {\Theta}^\infty_{k_0}, \\
y_k&=H_kx_k, \ k=0,1,2,\cdots
\end{split}
\right. \label{eq sssscdf}
\end{equation}
is ESMS, so $(F_k,G_k|H_k)$ is exactly detectable.
\end{IEEEproof}

\begin{remsub}\label{rem:3.1.2}\rm\
When system (\ref{eq vcgtgh1}) reduces to the deterministic time-invariant
system (\ref{eq vgthgh}), the uniform detectability, ${\cal K}^{n-1}$-exact
detectability and ${\cal K}^{\infty}$-exact detectability coincide with the
classical detectability of linear systems \cite{c1ss}.
Examples~\ref{exm:3.1.4}--\ref{exm:3.1.5} show that there is no inclusion
relation between uniform detectability and ${\cal K}^N$-exact
detectability for some $N>0$. We conjecture that if $(F_k,G_k|H_k)$ is
uniformly detectable, then there is a sufficiently large $N>0$ such that
$(F_k,G_k|H_k)$ is ${\cal K}^N$-exactly detectable.
\end{remsub}

Corresponding to Theorem~\ref{thm:2.1}, we also have the following theorem
for exact detectability, but its proof is very simple.

\begin{thmsub}\label{thm:3.1.1}\
If $(F_k,G_k|H_k)$ is ${\cal K}^{\infty}$-exactly detectable, then so is
$(F_k+M_kK_kH_k, G_k+N_kK_kH_k|H_k)$ for any output feedback $u_k=K_ky_k$.
\end{thmsub}

\begin{IEEEproof}\
We prove this theorem by contradiction. Assume that
$(F_k+M_kK_kH_k, G_k+N_kK_kH_k|H_k)$ is not ${\cal K}^{\infty}$-exactly
detectable. By Definition~\ref{def:3.1.3}, for system (\ref{eq vcgfdef}),
although the measurement equation becomes $y_k=H_kx_k\equiv 0$ for
$k\in {\cal N}_0$, the state equation
\begin{equation}
x_{k+1}=(F_k+M_kK_kH_k)x_k+(G_k+N_kK_kH_k)x_kw_k \label{eq cvghvh}
\end{equation}
is not ESMS. In view of $y_k=H_kx_k\equiv 0$, (\ref{eq cvghvh}) is
equivalent to
\begin{equation}
x_{k+1}=F_kx_k+G_kx_kw_k. \label{eq dvbfhd}
\end{equation}
Hence, under the condition of $y_k=H_kx_k\equiv 0$ for $k=0,1, 2, \cdots$,
if (\ref{eq cvghvh}) is not ESMS, then so is (\ref{eq dvbfhd}), which
contradicts the ${\cal K}^{\infty}$-exact detectability of $(F_k,G_k|H_k)$.
\end{IEEEproof}

It should be pointed out that Theorem~\ref{thm:3.1.1} does not hold for
${\cal K}^{N}$-exact detectability. That is, even if $(F_k,G_k|H_k)$ is
${\cal K}^{N}$-exactly detectable for $N\ge 0$,
$(F_k+M_kK_kH_k, G_k+N_kK_kH_k|H_k)$ may not be so, and such a counterexample
can be easily constructed.

\begin{propsub}\label{prop:3.1.1}\
If there exists a matrix sequence $\{K_k, k=0, 1, \cdots, \}$ such that
\begin{equation}
x_{k+1}=(F_k+K_kH_k)x_k+G_kx_kw_k \label{eq dvbfhdvgbg}
\end{equation}
is ESMS, then $(F_k,G_k|H_k)$ is ${\cal K}^{\infty}$-exactly detectable.
\end{propsub}

\begin{IEEEproof}\
Because (\ref{eq dvbfhdvgbg}) is ESMS, by Proposition~\ref{prop:2.1}
and Lemma~\ref{lem:3.1.3}, $(F_k+K_kH_k, G_k|H_k)$ is ${\cal
K}^{\infty}$-exactly detectable. By Theorem~\ref{thm:3.1.1}, for any
matrix sequence $\{L_k, k=0, 1, \cdots, \}$, $(F_k+K_kH_k+L_kH_k,
G_k|H_k)$ is also ${\cal K}^{\infty}$-exactly detectable. Taking
$L_k=-K_k$, we obtain that $(F_k,G_k|H_k)$ is ${\cal
K}^{\infty}$-exactly detectable. Thus, this proposition is shown.
\end{IEEEproof}

\begin{remsub}\label{rem:3.1.3}\rm\
In some previous references such as \cite{c2ddss, Ungureanu1}, if
system (\ref{eq dvbfhdvgbg}) is ESMS for some matrix sequence
$\{K_k\}_{k\in {\cal N}_0}$, then $(F_k,G_k|H_k)$ is called
stochastically detectable or detectable in conditional mean
\cite{Ungureanu1}. Proposition~\ref{prop:3.1.1} tells us that
stochastic detectability implies ${\cal K}^{\infty}$-exact
detectability, but the converse is not true. Such a counterexample
can be easily constructed; see the following
Example~\ref{exm:3.1.6}. The ${\mathcal K}^{\infty}$-exact
detectability implies that any  $k_0^\infty$-unobservable initial
state $\xi$ leads to an exponentially stable trajectory for any
$k_0\ge 0$. However, in the time-invariant system (\ref{eq
vcgdsdsf}), the stochastic detectability of (\ref{eq vcgdsdsf}) (or
$(F,G|H)$ for short) is equivalent to that there is a constant
output feedback gain matrix $K$, rather than necessarily a
time-varying feedback gain matrix sequence $\{K_k\}_{k\in {\mathcal
N}_0}$, such that
\begin{equation}
x_{k+1}=(F+KH)x_k+Gx_kw_k \label{eq cggfl}
\end{equation}
is ESMS; see \cite{c2ddss}.
\end{remsub}

\begin{exmsub}\label{exm:3.1.6}\rm\
Let $G_k=3$ for $k\ge 0$, and
$$
F_k=H_{k}=\left\{
\begin{array}{l}
1, \ \mbox{for} \ k=3n, \ n=1, 2,\cdots,\\
0, \ \mbox{otherwise}.
 \end{array}
 \right.
$$
By Lemma~\ref{lem:2.2}, for any output feedback $u_k=K_ky_k$, we have
$Ex_k^2=3^{(k-k_0)}Ex_{k_0}^2$ for $k>k_0$, where $x_k$ is the closed-loop
trajectory of
$$
x_{k+1}=(F_k+K_kH_k)x_k+3x_kw_k ,
$$
which is not ESMS. So $(F_k, G_k|H_k)$ is not stochastically detectable.
However, $(F_k,G_k|H_k)$ is not only ${\cal K}^\infty$- but also
${\cal K}^3$-exactly detectable, and $0$ is the unique $k^3$-unobservable
state.
\end{exmsub}

\begin{remsub}\label{rem:3.1.4}\rm\
According to the linear system theory, for the deterministic linear
time-invariant system (\ref{eq vgthgh}), the ${\cal K}^{\infty}$- and
${\cal K}^{n-1}$-exact detectability are equivalent. By the
${\cal H}$-representation theory \cite{xsdsffdf}, for (\ref{eq vcgdsdsf}),
the ${\cal K}^{\infty}$- and ${\cal K}^{[\frac {n(n+1)}2-1]}$-exact
detectability are also equivalent. So, in what follows, system
(\ref{eq vcgdsdsf}) (or $(F,G|H)$) is simply called exactly detectable.
\end{remsub}

\begin{remsub}\label{rem:3.1.5}\rm\
In Example~\ref{exm:3.1.3}, $(F_k|H_k)$ is stochastically detectable, but it
is not ${\cal K}^N$-exactly detectable for any $N\ge 0$. In
Example~\ref{exm:3.1.6}, $(F_k|H_k)$ is not stochastically detectable, but
it is ${\cal K}^N$-exactly detectable for $N\ge 3$. Hence, it seems that there
is no inclusion relation between stochastic detectability and
${\cal K}^N$-exact detectability.
\end{remsub}

\subsection{Lyapunov-Type Theorems under Exact Detectability\label{sec:LTTED-3}}

At present, we do not know whether Theorem~\ref{thm:4.1.2} holds
under exact detectability, but we are able to prove a similar result
to Theorem~\ref{thm:4.1.2} for a periodic system, namely, in
(\ref{eq vcgtgh1}), $F_{k+\tau}=F_k$, $G_{k+\tau}=G_k$,
$H_{k+\tau}=H_k$. Periodic systems are a class of very important
time-varying systems, which have been studied by many researchers;
see \cite{bittanti, c2ddss, dragan2013}.

\begin{thmsub}[Lyapunov-Type Theorem]\label{thm:5.3.1}\
Assume that system (\ref{eq vcgtgh1}) is a periodic system with the
period $\tau>0$. If system (\ref{eq vcgtgh1}) is ${\cal
K}^{N}$-exactly detectable for any fixed $N\ge 0$ and
$\{P_k>0\}_{k\ge 0}$ is a positive definite matrix sequence which
solves GLE (\ref{eq vhvcxf}), then the periodic system (\ref{eq
vcgtgh1}) is ESMS.
\end{thmsub}

\begin{IEEEproof}\
By periodicity, $P_k=P_{k+\tau}$. Select an integer $\bar\kappa>0$
satisfying $\bar\kappa\tau-1\ge N$. For $\kappa\ge\bar\kappa$, we
introduce the following backward difference equation
\begin{equation}
\left\{
\begin{array}{l}
-P_0^{\kappa\tau-1}(k)+F_{k}^TP_0^{\kappa\tau-1}(k+1)F_k+G_k^TP_0^{\kappa\tau-1}(k+1)G_k+H_k^TH_k=0,\\
P_0^{\kappa\tau-1}({\kappa\tau})=0, \ k=0, 1, \cdots, \kappa\tau-1.
\end{array}
\right. \label{eq vffggjjj}
\end{equation}
Set $V_k=x_k^TP_kx_k$, then associated with (\ref{eq vffggjjj}), we
have
\begin{align}
EV_0-EV_{\kappa\tau} &=
x_0^TP_0x_0-E[x_{\kappa\tau}^TP_{\kappa\tau}x_{\kappa\tau}]
=x_0^TP_0x_0-E[x_{\kappa\tau}^TP_{0}x_{\kappa\tau}]\nonumber\\
&=
\sum_{i=0}^{\kappa\tau-1}E\|y_i\|^2=x_0^TP_0^{\kappa\tau-1}(0)x_0,
\label{eq cvgvgfgf}
\end{align}
where the last equality is derived by using the completing squares
technique. We assert that $P_0^{\kappa\tau-1}(0)>0$. Otherwise,
there exists a nonzero $x_0$ satisfying
$x_0^TP_0^{\kappa\tau-1}(0)x_0=0$ due to $P_0^{\kappa\tau-1}(0)\ge
0$. As so, by ${\cal K}^{N}$-exact detectability, (\ref{eq
cvgvgfgf}) leads to
\begin{align}
0=\sum_{i=0}^{\kappa\tau-1}E\|y_i\|^2
&\ge \lambda_{\min}(P_0)\|x_0\|^2-\lambda_{\max}(P_0)\beta\lambda^{\kappa\tau}\|x_0\|^2\nonumber\\
&=
(\lambda_{\min}(P_0)-\lambda_{\max}(P_0)\beta\lambda^{\kappa\tau})\|x_0\|^2
, \label{eq vgvhgf}
\end{align}
where $\beta>1$ and $0<\lambda<1$ are defined in (\ref{eq vcgv}). If
$\kappa$ is taken sufficiently large such that $\kappa\ge
\kappa_0>0$ with $\kappa_0>0$ being a minimal integer satisfying
$\lambda_{\min}(P_0)-\lambda_{\max}(P_0)\beta\lambda^{\kappa_0\tau}>0$,
then (\ref{eq vgvhgf}) yields $x_0=0$, which contradicts $x_0\ne 0$.

If we let $P_{(n-1)\kappa\tau}^{n\kappa\tau-1}((n-1)\kappa\tau+k)$
denote the solution of
$$
\left\{
\begin{array}{l}
-P_{(n-1)\kappa\tau}^{n\kappa\tau-1}((n-1)\kappa\tau+k)+F_{(n-1)\kappa\tau+k}^TP_{(n-1)\kappa\tau}^{n\kappa\tau-1}((n-1)\kappa\tau+k+1)F_{(n-1)\kappa\tau+k}\\
\ \ \ \ \ \ \ \ \ \
+G_{(n-1)\kappa\tau+k}^TP_{(n-1)\kappa\tau}^{n\kappa\tau-1}((n-1)\kappa\tau+k+1)G_{(n-1)\kappa\tau+k}
 +H_{(n-1)\kappa\tau+k}^TH_{(n-1)\kappa\tau+k}=0,\\
P_{(n-1)\kappa\tau}^{n\kappa\tau-1}(n\kappa\tau)=0, \ k=0, 1,
\cdots, \kappa\tau-1; \ n=1,2,\cdots ,
\end{array}
\right.
$$
then by periodicity,
$P_{0}^{\kappa\tau-1}(0)=P_{(n-1)\kappa\tau}^{n\kappa\tau-1}((n-1)\kappa\tau)>0$,
and
\begin{align*}
EV_{(n-1)\kappa\tau}-EV_{n\kappa\tau} &=
\sum_{i=(n-1)\kappa\tau}^{n\kappa\tau-1}
E\|y_i\|^2=E[x_{(n-1)\kappa\tau}^TP_{(n-1)\kappa\tau}^{n\kappa\tau-1}((n-1)\kappa\tau)x_{(n-1)\kappa\tau}]\\
&=
E[x_{(n-1)\kappa\tau}^TP_{0}^{\kappa\tau-1}(0)x_{(n-1)\kappa\tau}]
\ge \varrho_0 \|x_{(n-1)\kappa\tau}\|^2 ,
\end{align*}
where $\varrho_0=\lambda_{\min} (P_{0}^{\kappa\tau-1})>0$.
Generally, for $0\le s\le \kappa\tau-1$, we define
$P_{(n-1)\kappa\tau+s}^{n\kappa\tau+s-1}((n-1)\kappa\tau+s+k)$ as
the solution to
$$
\left\{
\begin{array}{l}
-P_{(n-1)\kappa\tau+s}^{n\kappa\tau+s-1}((n-1)\kappa\tau+s+k)+F_{(n-1)\kappa\tau+s+k}^TP_{(n-1)\kappa\tau+s}^{n\kappa\tau+s-1}((n-1)\kappa\tau+s+k+1)
F_{(n-1)\kappa\tau+s+k}\\
\ \ \
+G_{(n-1)\kappa\tau+s+k}^TP_{(n-1)\kappa\tau+s}^{n\kappa\tau+s-1}((n-1)\kappa\tau+s+k+1)G_{(n-1)\kappa\tau+s+k}
 +H_{(n-1)\kappa\tau+s+k}^TH_{(n-1)\kappa\tau+s+k}=0,\\
P_{(n-1)\kappa\tau+s}^{n\kappa\tau+s-1}(n\kappa\tau+s)=0, \ k=0, 1,
\cdots, \kappa\tau-1; \ n=1,2,\cdots.
\end{array}
\right.
$$
It can be shown that
$P_{(n-1)\kappa\tau+s}^{n\kappa\tau+s-1}((n-1)\kappa\tau+s)=P_{s}^{\kappa\tau+s-1}(s)>0$
and
$$
\sum_{i=(n-1)\kappa\tau+s}^{n\kappa\tau+s-1}
E\|y_i\|^2=E[x_{(n-1)\kappa\tau+s}^TP_{s}^{\kappa\tau+s-1}(s)x_{(n-1)\kappa\tau+s}],
$$
provided that we take $\kappa\ge \max_{0\le s\le
\kappa\tau-1}{\kappa_s}$, where $\kappa_s>0$ is the minimal integer
satisfying
$\lambda_{\min}(P_s)-\lambda_{\max}(P_s)\beta\lambda^{\kappa_s\tau}>0$.

Summarizing the above discussions, for any $k\ge 0$ and
$\hat\kappa>\max\{\bar\kappa, \max_{0\le s\le
\kappa\tau-1}{\kappa_s}\}$, we have
$$
EV_k-EV_{k+{\hat \kappa}\tau}=\sum_{i=k}^{k+{\hat
\kappa}\tau-1}E\|y_i\|^2\ge \rho E\|x_k\|^2 ,
$$
where $\rho=\min_{0\le s\le {\hat \kappa}\tau-1} {\rho_s}>0$ with
$\rho_s=\lambda_{\min} [P_{s}^{{\hat\kappa}\tau+s-1}(s)]$. The rest
is similar to the proof of Theorem~\ref{thm:4.1.2} and thus is
omitted.
\end{IEEEproof}

In Theorem~\ref{thm:5.3.1}, if $\{P_k>0\}_{k\ge 0}$ is weaken as
$\{P_k\ge 0\}_{k\ge 0}$, then we have

\begin{thmsub}[Lyapunov-Type Theorem]\label{thm:5.3.2}\
Assume that system (\ref{eq vcgtgh1}) is a periodic system with the
period $\tau>0$. If (i) system (\ref{eq vcgtgh1}) is ${\cal
K}^{N}$-exactly detectable for any fixed $N\ge 0$; (ii)
$\{P_k\ge0\}_{k\ge 0}$ is a positive semi-definite matrix sequence
which solves GLE (\ref{eq vhvcxf}); (iii)
$\mbox{Ker}(P_0)=\mbox{Ker}(P_1)=\cdots=\mbox{Ker}(P_{\tau-1})$,
then the periodic system (\ref{eq vcgtgh1}) is ESMS.
\end{thmsub}

\begin{IEEEproof}\
From GLE (\ref{eq vhvcxf}), it is easy to show (e.g., see
Theorem~3.2 in \cite{zhang1}) that $\mbox{Ker}(P_k)\subset
\mbox{Ker}(H_k)$, $F_k\mbox{Ker}(P_k)\subset \mbox{Ker}(P_{k+1})$,
$G_k\mbox{Ker}(P_k)\subset \mbox{Ker}(P_{k+1})$. In addition, in
view of $\mbox{Ker}(P_0)=\cdots=\mbox{Ker}(P_{\tau-1})$ and
$P_{\tau+k}=P_k$, there is a common orthogonal matrix $S$ such that
for any $k\ge 0$, there hold
$$
S^TP_kS=\left[
\begin{array}{cc}
0 &  0\\
0 &  P^{22}_k
\end{array}
\right], \   P^{22}_k\ge 0,  \  S^TH_k^TH_kS=\left[
\begin{array}{cc}
0 &  0\\
0 & (H^{22}_k)^T H^{22}_k
\end{array}
\right],
$$
$$
  \  S^TF_kS=\left[
\begin{array}{cc}
F^{11}_k &  F^{12}_k\\
0 &  F^{22}_{k}
\end{array}
\right], \ \ S^TG_kS=\left[
\begin{array}{cc}
G^{11}_k &  G^{12}_k\\
0 &  G^{22}_k
\end{array}
\right].
$$
Pre- and post-multiplying $S^T$ and $S$ on both sides of GLE
(\ref{eq vhvcxf}) gives rise to
$$
-S^TP_kS+ S^TF^T_k S \cdot S^TP_{k+1}S\cdot S^TF_kS+ S^TG_k^TS\cdot
S^TP_{k+1}S \cdot S^TG_kS+S^TH_k^TS\cdot S^TH_kS=0 ,
$$
which is equivalent to
\begin{equation}
-P^{22}_k+(F^{22}_k)^TP^{22}_{k+1}F^{22}_k+(G^{22}_k)^TP^{22}_{k+1}G^{22}_k+(H^{22}_k)^TH^{22}_k=0.
\label{eq cvgvgh}
\end{equation}
Set $\eta_k=\left[
\begin{array}{c}
\eta_{1,k}\\
\eta_{2,k}
\end{array}
\right]=S^Tx_k= \left[
\begin{array}{cc}
S_{11} &  S_{12}\\
S_{21} &  S_{22}
\end{array}
\right]^Tx_k$, then it follows that
\begin{equation}
\left\{
\begin{split}
\eta_{1,k+1}&=F^{11}_k\eta_{1,k}+G^{11}_k\eta_{1,k}w_k+F^{12}_k\eta_{2,k}+G^{12}_k\eta_{2,k}w_k, \\
\eta_{2,k+1}&=F^{22}_k\eta_{2,k}+G^{22}_k\eta_{2,k}w_k, \\
y_k&=H_kS\eta_k.
\end{split}
\right. \label{eq sdfthtyhd}
\end{equation}
It can be easily seen that $y_k=H_kS\eta_k\equiv 0, a.s.$ iff
$H^{22}_k\eta_{2,k}\equiv 0, a.s.$, for which a sufficient condition
is $\eta_{2,k}=0$. By ${\cal K}^{N}$-exact detectability,
$\eta_{1,k+1}=F^{11}_k\eta_{1,k}+G^{11}_k\eta_{1,k}w_k$ is ESMS. To
show that $\eta_{2,k+1}=F^{22}_k\eta_{2,k}+G^{22}_k\eta_{2,k}w_k$ is
ESMS, we consider the following reduced-order state-measurement
equation
\begin{equation}
\left\{
\begin{split}
\eta_{2,k+1}&=F^{22}_k\eta_{2,k}+G^{22}_k\eta_{2,k}w_k , \\
{\bar y}_k&=H^{22}_k\eta_{2,k} .
\end{split}
\right. \label{eq vhvhhss}
\end{equation}
Obviously, system (\ref{eq vhvhhss}) is still a periodic system and
has the same period $\tau>0$ as (\ref{eq vcgtgh1}).

In the following, we show that (\ref{eq vhvhhss}) is also ${\cal
K}^N$-exactly detectable. Because system (\ref{eq vcgtgh1}) is
${\cal K}^{N}$-exactly detectable, for any $k\ge 0$, $y_i\equiv 0$
a.s. for $i=k, \cdots, k+N$, implies that there are constants
$\beta_0>1$ and $0<\lambda_0<1$ such that
\begin{equation}
E\|x_{t}\|^2=E\|\phi_{k+t,k}x_k\|^2\le \beta_0 E\|x_{k}\|^2
\lambda_0^{(t-k)},\ t\ge k \label{eq vcgwewe}
\end{equation}
for any $k^N$-unobservable state $x_k$. Take $x_k=S\eta_k=S\left[
\begin{array}{c}
0\\
\eta_{2,k}
\end{array}
\right], $ with $\eta_{2,k}$ being a $k^N$-unobservable state of
(\ref{eq vhvhhss}), i.e., ${\bar y}_i=H^{22}_k\eta_{2,i}=0$ for
$i=k,\cdots, k+N$. Then ${y}_i=H_iS\left[
\begin{array}{c}
0\\
\eta_{2,k}
\end{array}
\right]=0$ for $i=k, \cdots, k+N$. Hence, (\ref{eq vcgwewe}) holds.
Substituting $x_k=S\left[
\begin{array}{c}
0\\
\eta_{2,k}
\end{array}
\right] $ into (\ref{eq vcgwewe}) yields
\begin{equation}
E\|\eta_{2,t}\|^2\le \beta_0 E\|\eta_{2,k}\|^2 \lambda_0^{(t-k)}, \
t\ge k.  \label{eq vcgweweb}
\end{equation}
So (\ref{eq vhvhhss}) is ${\cal K}^{N}$-exactly detectable.

Associated with (\ref{eq vhvhhss}), the GLE (\ref{eq cvgvgh}) admits
a positive definite solution sequence $\{P_k>0\}_{k\ge0}$. Applying
Theorem~\ref{thm:5.3.1}, the subsystem (\ref{eq vhvhhss}) is ESMS.
Since $ \eta_{1,k+1}=F^{11}_k\eta_{1,k}+G^{11}_k\eta_{1,k}w_k $ has
been shown to be ESMS, there are constants $\beta_1>1$ and
$0<\lambda_1<1$ such that
\begin{equation}
E\|\eta_{1,t}\|^2\le \beta_1 E\|\eta_{1,k}\|^2 \lambda_1^{(t-k)}, \
t\ge k.  \label{eq vcgwhhhdv}
\end{equation}
Set $\beta:=\max\{\beta_0, \beta_1\}$,
$\lambda=\max\{\lambda_0,\lambda_1\}$, then the composite system
(\ref{eq sdfthtyhd}) is ESMS with
$$
E\|\eta_{t}\|^2=E\|\eta_{1,t}\|^2+E\|\eta_{2,t}\|^2\le \beta
E\|\eta_{k}\|^2 \lambda^{(t-k)}, \ t\ge k ,
$$
which deduces that the periodic system (\ref{eq vcgtgh1}) is ESMS
because (\ref{eq vcgtgh1}) and (\ref{eq sdfthtyhd}) are equivalent.
\end{IEEEproof}

Finally, we consider the linear time-invariant stochastic system
(\ref{eq vcgdsdsf}) and present a Lyapunov-type theorem as a
complementary result of Theorem 19 \cite{lizao}. Associated with
(\ref{eq vcgdsdsf}), we introduce the linear symmetric operator
${\cal L}_{F,G}$, called the generalized Lyapunov operator (GLO), as
follows:
$$
{\cal L}_{F,G} Z=FZF^T+GZG^T, \ \ Z\in {\cal S}_n.
$$
Moreover, for system (\ref{eq vcgdsdsf}), the GLE (\ref{eq vhvcxf}) becomes
\begin{equation}
-P+F^TPF+G^TPG+H^TH=0. \label{eq vhdfvf}
\end{equation}

\begin{thmsub}\label{thm:5.1.2}\
Suppose that $\sigma ({\cal L}_{F,G})\subset \bar {\odot}
:=\{\lambda: |\lambda|\le 1 \}$ and $(F,G|H)$ is exactly detectable.
If $P$ is a real symmetric solution of (\ref{eq vhdfvf}), then $P\ge
0$ and $(F, G)$ is stable, i.e., the state trajectory of (\ref{eq
vcgdsdsf}) is asymptotically mean square stable.
\end{thmsub}

In order to prove Theorem~\ref{thm:5.1.2}, we need to cite the
well-known Krein-Rutman Theorem as follows:

\begin{lemsub}[see \cite{Schneider}]\label{lem:5.1.1}\
Let $\beta:=\max_{\lambda_i\in \sigma ({\cal L}_{F,G})} |\lambda_i|$
be the spectral radius of ${\cal L}_{F,G}$. Then there exists a
nonzero $X\ge 0$ such that ${\cal L}_{F,G}X=\beta X$.
\end{lemsub}

\begin{IEEEproof}[Proof of Theorem~\ref{thm:5.1.2}]\
Because $\sigma ({\cal L}_{F,G})\subset \bar {\odot}$, the spectral
radius $\beta\le 1$. If $\beta<1$, then this means that $(F, G)$ is
stable by \cite[Lemma 3]{lizao}, which yields $P\ge 0$ according
to \cite[Lemma 17]{lizao}. If $\beta=1$, then by Lemma~\ref{lem:5.1.1},
there exists a nonzero $X\ge 0$, such that ${\cal L}_{F,G}X=X$.
Therefore, we have
\begin{align}
0\ge \langle -H^TH, X\rangle &= \langle -P+{\cal L}^*_{F,G}(P),
X\rangle
= \langle -P,X\rangle + \langle P, {\cal L}_{F,G}(X)\rangle \nonumber\\
&= \langle -P,X\rangle + \langle P,X\rangle = \langle 0, X\rangle =0
, \label{eq vhgvh}
\end{align}
where $\langle A,B\rangle :=\mathrm{trace}(A^TB)$, ${\cal
L}^*_{F,G}$ is the adjoint operator of ${\cal L}_{F,G}$, and ${\cal
L}^*_{F,G}(P)=F^TPF+G^TPG$. From (\ref{eq vhgvh}) it follows that
$\mathrm{trace}(H^THX)=0$, which implies $HX=0$ due to $X\ge 0$.
However, according to \cite[Theorem 8-(4)]{lizao}, ${\cal
L}_{F,G}X=X$ together with $HX=0$, contradicts the exact
detectability of $(F,G|H)$. Hence, we must have $0\le \beta<1$, and
this theorem is verified.
\end{IEEEproof}

\begin{remsub}\label{rem:5.1.1}\rm\
Following the line of Theorem~\ref{thm:5.1.2}, Conjecture 3.1 in
\cite{xshf} can also be verified.
\end{remsub}

%%%%%%%%%%%%%%%%%%%%%%%%%%%%%%%%%%%%%%%%%%%%%%%%%%%%%%%%%%%%%%%%%%%%%%%%%%%%%%%%
\section{Exact Observability \label{sec:TLTT}}

This section introduce another definition called ``exact
observability'' for system (\ref{eq vcgtgh1}), which is stronger
than exact detectability and  also coincides with the classical
observability when system (\ref{eq vcgtgh1}) reduces to the
deterministic linear time-invariant system (\ref{eq vgthgh}).

We first give the following definitions:

\begin{definition}\label{def:3.2.1}\
System (\ref{eq vcgtgh1}) is called $k_0^\infty$-exactly observable if
$x_{k_0}=0$ is the unique $k_0^\infty$-unobservable state. Similarly,
system (\ref{eq vcgtgh1}) is called $k_0^{s_0}$-exactly observable if
$x_{k_0}=0$ is the unique $k_0^{s_0}$-unobservable state.
\end{definition}

\begin{definition}\label{def:3.2.2}\
System (\ref{eq vcgtgh1}) (or $(F_k,G_k|H_k)$) is said to be
${\cal K}^\infty$-exactly observable if for any $k\in [0,\infty)$,
system (\ref{eq vcgtgh1}) is $k^\infty$-exactly observable. If for any
time $k\in [0,\infty)$, there exists a nonnegative integer $N\ge 0$ such
that system (\ref{eq vcgtgh1}) is $k^{N}$-exactly observable, then system
(\ref{eq vcgtgh1}) (or $(F_k,G_k|H_k)$) is said to be ${\cal K}^{N}$-exactly
observable. Similarly, the ${\cal K}^{WFT}$- and ${\cal K}^{FT}$-exact
observability can be defined.
\end{definition}

Combining Lemmas~\ref{lem:3.1.1}--\ref{lem:3.1.2} together, a sufficient
condition for the exact observability is presented as follows.

\begin{theorem}\label{thm:3.2.1}\
If ${\mbox {rank}}H_{\infty,k}=n$ for any $k\ge 0$, then
$(F_k,G_k|H_k)$ is ${\cal K}^\infty$-exactly observable. In
particular, if ${\mbox {rank}}H_{k+s_0,k}=n$ for some fixed integer
$s_0\ge 0$ and any $k\ge 0$, then system (\ref{eq vcgtgh1}) is not
only ${\cal K}^\infty$- but also ${\cal K}^{s_0}$-exactly
observable. Here $H_{l,k}$ is defined in Lemma~\ref{lem:2.1}, and
$$
H_{\infty,k}=\left[
\begin{array}{cc}
H_k\\
(I_2\otimes H_{k+1})\phi_{k+1,k}\\(I_{2^2}\otimes
H_{k+2})\phi_{k+2,k}\\ \vdots \\ (I_{2^{l-k}}\otimes
H_l)\phi_{l,k}\\ \vdots
\end{array}
\right]. \label{eq 2bv}
$$
\end{theorem}

The next corollary follows immediately from Theorem~\ref{thm:3.2.1}.

\begin{corollary}\label{cor:3.2.1}\
If $H_{k}$ is nonsingular for $k\ge 0$, then system (\ref{eq vcgtgh1}) is
${\cal K}^{0}$-exactly observable.
\end{corollary}

By Definitions~\ref{def:3.2.1}--\ref{def:3.2.2}, $k_0^\infty$
(resp. $k_0^{s_0}$)-exact observability is stronger than $k_0^\infty$
(resp. $k_0^{s_0}$)-exact detectability. Likewise, ${\cal K}^\infty$
(resp. ${\cal K}^{WFT}$, ${\cal K}^{FT}$, ${\cal K}^{N}$)-exact observability
is stronger than ${\cal K}^\infty$ (resp. ${\cal K}^{WFT}$, ${\cal
K}^{FT}$, ${\cal K}^{N}$)-exact detectability. A necessary and sufficient
condition for the ${\cal K}^{N}$-exact observability was presented in
\cite{xsdsffdf} based on the $\cal H$-representation theory developed therein.
Below, we give another equivalent theorem based on Lemma~\ref{lem:2.1}.

\begin{theorem}\label{thm:3.2.2}\
(i) System (\ref{eq vcgtgh1}) is ${\cal K}^{N}$-exactly observable
iff for any $k\in {\cal N}_0$, the Gramian ${\cal O}_{k+N,k}$ is a
positive definite matrix. (ii) If system (\ref{eq vcgtgh1}) is
${\cal K}^{WFT}$-exactly observable and $\{P_k\ge 0\}_{k\ge 0}$
solves the GLE (\ref{eq vhvcxf}), then $P_k>0$ for any $k\ge 0$.
\end{theorem}

\begin{IEEEproof}\
We note that $y_i\equiv 0$ a.s. for $i=k, k+1, \cdots, k+N$, is
equivalent to $\sum_{i=k}^{k+N} E\|y_i\|^2=0$. By
Lemma~\ref{lem:2.1}, $\sum_{i=k}^{k+N} E\|y_i\|^2=E[x_k^T{\cal
O}_{k+N,k}x_k]=0$. So system (\ref{eq vcgtgh1}) is exactly
observable iff $\sum_{i=k}^{k+N} E\|y_i\|^2=E[x_k^T{\cal
O}_{k+N,k}x_k]=0$ implies $x_k=0$ a.s., which is equivalent to
${\cal O}_{k+N,k}>0$ due to ${\cal O}_{k+N,k}\ge 0$. Hence, (i) is
proved.

Now we prove (ii) by contradiction. If some $P_{k_0}$ is not
strictly positive definite, then there exists a nonzero $x_{k_0}\in
l^2_{{\cal F}_{k_0-1}}$ such that $E[x_{k_0}^TP_{k_0}x_{k_0}]=0$. By
the ${\cal K}^{WFT}$-exact observability of $(F_k,G_k|H_k)$, there
is $s_0\ge 0$ such that system (\ref{eq vcgtgh1}) is
$k_0^{s_0}$-exactly observable. Since the following identity
\begin{eqnarray}
E[x_k^TP_kx_k]-E[x_{s+1}^TP_{s+1}x_{s+1}]&=&
\sum_{i=k}^{s}E\|y_i\|^2 \label{eq fsdfddd}
\end{eqnarray}
holds for any $s\ge k\ge 0$, it follows that
$$
0\le
\sum_{i=k_0}^{k_0+s_0}E\|y_i\|^2=-E[x_{k_0+s_0+1}^TP_{k_0+s_0+1}x_{k_0+s_0+1}]\le
0
$$
and accordingly $y_i\equiv 0$ a.s. for $i\in [k_0,k_0+s_0]$. By the
$k_0^{s_0}$-exact observability, $x_{k_0}=0$, which contradicts
$x_{k_0}\ne 0$. Hence, (ii) is proved.
\end{IEEEproof}

\begin{remark}\label{rem:3.2.1}\rm\
Theorem~\ref{thm:3.2.2}-(i) shows that the ${\cal K}^{N}$-exact observability
is weaker than the uniform observability given in \cite{c2ddss}, where it was
proved that system (\ref{eq vcgtgh1}) is uniformly observable iff there are
$N\ge 0$ and $\gamma>0$ such that ${\cal O}_{k+N,k}\ge \gamma I$ for any
$k\in {\cal N}_0$.
\end{remark}

\begin{remark}\label{rem:3.2.2}\rm\
There is no inclusion relation between uniform detectability and exact
observability. For example, in Example~\ref{exm:3.1.4}, $(F_k|H_k)$ is
${\cal K}^N$-exactly observable, but it is not uniformly detectable. On the
other hand, in Example~\ref{exm:3.1.5}, $(F_k|H_k)$ is uniformly detectable,
but it is not ${\cal K}^0$-exactly observable.
\end{remark}

Similar to exact detectability, we also have the following inclusion
relation for exact observability:
\begin{align*}
& {\cal K}^{N}\mbox {-exact observability} ~\Longrightarrow~
{\cal K}^{FT}\mbox{-exact observability} \\
\Longrightarrow~& {\cal K}^{WFT}\mbox{-exact observability}
~\Longrightarrow~ {\cal K}^{\infty}\mbox{-exact observability}.
\end{align*}

The following Lyapunov-type theorem can be viewed as a corollary of
Theorem~\ref{thm:5.3.1}.

\begin{theorem}[Lyapunov-Type Theorem]\label{thm:4.2.1}\
Assume that system (\ref{eq vcgtgh1}) is a periodic system with the
period $\tau>0$. If system (\ref{eq vcgtgh1}) is ${\cal
K}^{N}$-exactly observable for any $N\ge 0$ and $\{P_k\ge 0\}_{k\ge
0}$ solves GLE (\ref{eq vhvcxf}), then the periodic system (\ref{eq
vcgtgh1}) is ESMS.
\end{theorem}

\begin{IEEEproof}\
By Theorem~\ref{thm:3.2.2}-(ii), $P_k>0$ for $k\ge 0$. Because
${\cal K}^{N}$-exact observability must be ${\cal K}^{N}$-exact
detectability, this theorem is an immediate corollary of
Theorem~\ref{thm:5.3.1}.
\end{IEEEproof}

%%%%%%%%%%%%%%%%%%%%%%%%%%%%%%%%%%%%%%%%%%%%%%%%%%%%%%%%%%%%%%%%%%%%%%%%%%%%%%%%
\section{Additional Comments\label{sec:LTTED}}
%%%%%%%%%%%%%%%%%%%%%%%%%%%%%%%%%%%%%%%%%%%%%%%%%%%%%%%%%%%%%%%%%%%%%%%%%%%%%%%%

At the end of this paper, we give the following comments:

\begin{itemize}

\item[(i)]
In \cite{ni, zhangtan, shenl}, exact observability and
detectability of linear stochastic time-invariant systems with
Markov jump were studied. How to extend various definitions of
this paper to linear time-varying Markov jump systems is an
interesting research topic that merits further study.

\item[(ii)]
Following the line of \cite{xsdsffdf} that transforms the system
(\ref{eq vcgtgh1}) into a deterministic time-varying system, it is
easy to give some testing criteria for uniform detectability and
observability of (\ref{eq vcgtgh1}) by means of the existing results on
deterministic time-varying systems\cite{peters}. In addition,
applying the infinite-dimensional operator theory, the spectral
criterion for stability of system (\ref{eq vcgtgh1}) is also a
valuable research issue.

\item[(iii)]
In view of Remarks~\ref{rem:3.1.3}--\ref{rem:3.1.4}, we know that,
for linear time-invariant system (\ref{eq vcgdsdsf}), stochastic
detectability implies exact detectability. In \cite{zhangtan}, it
was shown that exact detectability is equivalent to the so-called
``W-detectability" (see \cite[Definition 3]{zhangtan}). A new
definition called ``weak detectability" was introduced in
\cite{Ungureanu2}, where a counter-example (see Example 15 in
\cite{Ungureanu2}) shows that W-detectability does not imply weak
detectability. In particular, it was proved in \cite{Ungureanu2} that
weak detectability can be derived from stochastic detectability. It
is easy to prove that weak detectability implies exact
detectability. In summary, we have the following inclusion relation:
\begin{align*}
\mbox{stochastic detectability} \Rightarrow \mbox{weak
detectability}\Rightarrow \mbox{exact detectability}\Leftrightarrow
\mbox{W-detectability}.
\end{align*}
As stated in \cite{Ungureanu2}, the converse implication that
whether W-detectability or exact detectability implies weak
detectability is an open question.

\item[(iv)] Lemmas~\ref{lem:2.2}--\ref{lem:2.1} are important,
which will have potential applications in mean stability analysis
and system synthesis.

\item[(v)] This paper reveals some essential differences between
linear time-varying and time-invariant systems. For example, for
linear time-invariant system (\ref{eq vcgdsdsf}), exact
detectability and exact observability can be uniquely defined, but
they exhibit diversity for LDTV system (\ref{eq vcgtgh1}).
Moreover, many equivalent relations in linear time-invariant system
(\ref{eq vcgdsdsf}) such as
\begin{align*}
\mbox{uniform detectability} \Leftrightarrow \mbox{exact
detectability}, \qquad \mbox{uniform observability} \Leftrightarrow
\mbox{exact observability}
\end{align*}
do not hold for LDTV system (\ref{eq vcgtgh1}).

\end{itemize}

%%%%%%%%%%%%%%%%%%%%%%%%%%%%%%%%%%%%%%%%%%%%%%%%%%%%%%%%%%%%%%%%%%%%%%%%%%%%%%%%
\section{Conclusion\label{sec:Concl}}
%%%%%%%%%%%%%%%%%%%%%%%%%%%%%%%%%%%%%%%%%%%%%%%%%%%%%%%%%%%%%%%%%%%%%%%%%%%%%%%%

This paper has introduced the new concepts on detectability and
observability for LDTV stochastic systems with multiplicative noise.
Uniform detectability defined in this paper can be viewed as an
extended version of that in \cite{bdo}. Various definitions on exact
detectability and observability are extensions of those in
\cite{Damm1,xshf, lizao, zhang1, c3dd} to LDTV stochastic systems.
Different from time-invariant systems, defining exact detectability
and exact observability for the time-varying stochastic system
(\ref{eq vcgtgh1}) is much more complicated. We have also  obtained
some Lyapunov-type theorems under uniform detectability of LDTV
systems, ${\cal K}^N$-exact detectability and ${\cal K}^N$-exact
observability of linear discrete periodic systems. We believe that
all these new concepts that have been introduced herein will play
important roles in control and filtering design of LDTV systems.

%%%%%%%%%%%%%%%%%%%%%%%%%%%%%%%%%%%%%%%%%%%%%%%%%%%%%%%%%%%%%%%%%%%%%%%%%%%%%%%%
\section*{Acknowledgement}
%%%%%%%%%%%%%%%%%%%%%%%%%%%%%%%%%%%%%%%%%%%%%%%%%%%%%%%%%%%%%%%%%%%%%%%%%%%%%%%%

The first author would be grateful to his
doctoral students Yong Zhao and Yaning Lin for their valuable
discussions in proving Lemmas~\ref{lem:2.2}--\ref{lem:2.1}.

%%%%%%%%%%%%%%%%%%%%%%%%%%%%%%%%%%%%%%%%%%%%%%%%%%%%%%%%%%%%%%%%%%%%%%%%%%%%%%%%
%%%%%%%%%%%%%%%%%%%%%%%%%%%%%%%%%%%%%%%%%%%%%%%%%%%%%%%%%%%%%%%%%%%%%%%%%%%%%%%%

\end{document}